\documentclass[12pt,a4paper,fleqn]{article}
\usepackage{a4wide,amsfonts,amsmath,latexsym,amssymb,euscript,graphicx,units,mathrsfs}

\usepackage{graphicx}
\usepackage{color}
\usepackage{amssymb}
\usepackage{amssymb}
\usepackage[T1]{fontenc}
\usepackage{latexsym}
\usepackage{xypic}
\usepackage{eufrak}
\usepackage{euscript}
\usepackage{amsfonts,amsmath}
\usepackage{verbatim}
\usepackage{fancyhdr}
\usepackage[english]{babel}
\usepackage{mathrsfs}
\usepackage{units}
\usepackage{hyperref}
\usepackage{tikz}
\usepackage{subfloat}
\usepackage[lofdepth,lotdepth]{subfig}

\newtheorem{prop}{Proposition}[section]
\newtheorem{cor}[prop]{Corollary}
\newtheorem{lemme}[prop]{Lemma}
\newtheorem{rem}[prop]{Remark}

\newtheorem{thm}[prop]{Theorem}
\newtheorem{defi}[prop]{Definition}
\newtheorem{example}[prop]{Example}

\newtheorem{as}[prop]{Assumption}
\newtheorem{prob}[prop]{Problem}

\renewcommand{\geq}{\geqslant}
\def\leq{\leqslant}

\newcommand{\R}{\mathbb{R}}

\def\var{{\mathbb{Var}}}

\def\1{{\mathbf{1}}}

\def\1{{\mathbf{1}}}
\def\0.5{{\frac{1}{2}}}
\def\var{{\rm{Var}}}

\newcommand{\fin}
{ \vspace{-0.6cm}
\begin{flushright}
\mbox{$\Box$}
\end{flushright}
\noindent }

\newcommand{\qed}{\nopagebreak\hspace*{\fill}
{\vrule width6pt height6ptdepth0pt}\par}

\usepackage{fullpage}

\usepackage{hyperref}

\newcommand{\leb}{\ell}

\renewcommand{\phi}{\varphi}
\renewcommand{\kappa}{\varkappa}

\newcommand{\fact }{\overline}

\newlength{\querylen}
\usepackage{fancybox}

\begin{document}

\begin{center}
{\Large{\bf Fine Gaussian fluctuations on the Poisson space, I:\\ {contractions, cumulants and geometric random graphs}}}
\normalsize
\\~\\ by Rapha\"el Lachi\`eze-Rey and Giovanni Peccati \\ {\it  Universit\'e Paris Descartes et Universit\'e du Luxembourg}\\~\\
\end{center}

{\small \noindent {\bf Abstract}: We study the normal approximation of functionals of Poisson measures having the form of a finite sum of multiple integrals. When the integrands are nonnegative, our results yield necessary and sufficient conditions for central limit theorems. These conditions can always be expressed in terms of contraction operators or, equivalently, fourth cumulants. Our findings are specifically tailored to deal with the normal approximation of the geometric $U$-statistics introduced by Reitzner and Schulte (2011). In particular, we shall provide a new analytic characterization of geometric random graphs whose edge-counting statistics exhibit asymptotic Gaussian fluctuations, and describe a new form of Poisson convergence for stationary random graphs with sparse connections. In a companion paper, the above analysis is extended to general $U$-statistics of marked point processes with possibly rescaled kernels.\\



\section{Introduction}

This paper concerns the normal approximation of random variables living inside a fixed sum of Wiener chaoses associated with a Poisson measure over a Borel measure space. Our main theoretical tools come from the two papers \cite{PSTU, PecZheng}, respectively by Peccati {\it et al.} and Peccati and Zheng, where the normal approximation of functional of Poisson measures is studied by combining two probabilistic techniques, namely the {\it Stein's method}  and the {\it Malliavin calculus of variations}. 

We shall focus on conditions implying that a given sequence of random variables satisfies a central limit theorem (CLT), where the convergence in distribution takes place in the sense of the Wasserstein distance (see Section \ref{ss:intropoisson} for definitions). Our main concern is to provide {\it analytic conditions} for asymptotic normality, that is, conditions only involving expressions related to the kernels in the chaotic expansion of a given random variable. In particular, our approach does not involve computations based on the method of moments and cumulants (with the exception of Theorem \ref{t:poisson}, where we deal with Poisson approximations). 

\medskip

The main contributions of our paper are the following:

\begin{itemize}

\item[--] In Theorem \ref{t:main1}, we shall prove that conditions for asymptotic normality can be expressed in terms of norms of {\it contraction operators} (see Section \ref{ss:star}). These analytic objects already appear in CLTs living inside a fixed Wiener chaos (see \cite{PSTU, PecZheng}), and are a crucial tool in order to effectively assess bounds based on Malliavin operators. One further important point is that the use of contraction operators allows one to neatly distinguish the contribution of each chaotic projection to the CLT, as well as to deduce joint CLTs for these projections starting from the asymptotic normality of their sum (see Proposition \ref{p:sumtosingle}). 

\item[--] In Theorem \ref{t:main2} we shall prove that, when specialized to random variables such that each kernel in the Wiener-It\^o representation has a constant sign, our results yield {\it necessary and sufficient conditions} for asymptotic normality.  The main tools in order to show such a result are two new analytic bounds, stated in Proposition \ref{p:1bnd} and Proposition \ref{p:2bnd}. These findings extend to the Poisson framework the `fourth moment theorem' proved by Nualart and Peccati (in a Gaussian setting) in \cite{nuapec}, a result that has been the starting point of a new line of research in stochastic analysis -- see the book \cite{np-book}, as well as the constantly updated webpage $${\tt http://www.iecn.u-nancy.fr/~nourdin/steinmalliavin.htm.}$$ 

\item[--] As discussed below, random variables having Wiener-It\^o kernels with constant sign appear quite naturally in problems arising in stochastic geometry. In particular, we shall use our results in order to provide an exhaustive characterization of stationary geometric random graphs whose edge counting statistics exhibit asymptotic Gaussian fluctuations (see Theorem \ref{thm:CLT-stat-graphs}). This family of geometric graphs contains e.g. {\it interval graphs} and {\it disk graphs} -- see e.g.  \cite{CCJ, DFR, FV, McDiarmid, Muller, penrosebook}.  Our characterization of geometric random graphs involves `diagonal subsets' of Cartesian products, that are reminiscent of the combinatorial conditions for CLTs used by Blei and Janson in \cite{bleijanson}, in the framework of CLTs for finite Rademacher sums (see also \cite[Section 6]{npr_ejp}). As a by-product of our analysis (see Theorem \ref{t:poisson}), we shall illustrate a new form of Poisson convergence for random graphs with sparse connections.

\end{itemize}

We stress that one of our main motivations comes from a remarkable paper by Reitzner and Schulte \cite{lesmathias}, laying the foundations of a general theory for CLTs involving $U$-statistics based on Poisson point processes. In particular, one of the crucial insights of  \cite{lesmathias} concerns the use of a formula by Last and Penrose (see \cite{lastpenrose}), providing explicit expressions for Wiener-It\^o chaotic decompositions in terms of difference operators (see Theorem \ref{t:ustat}). It is interesting to note that Last and Penrose's formula is the Poisson analogous of the so-called `Stroock formula' of Malliavin calculus, which is in turn an important tool  for proving CLTs involving  non-linear functionals of Gaussian measures (see e.g. \cite[Corollary 2.7.8]{np-book} for a discussion of this point). We shall see that our findings complement and extend the results proved in \cite{lesmathias} in several directions. See also Decreusefond {\it et al.} \cite{DFR}, Ferraz and Vergne \cite{FV}, Last {\it et al.} \cite{lpst}, Minh \cite{minh}, Schulte \cite{schulte2011}, Schulte and Thaele \cite{SchTh, SchTh2012},  for several new findings pertaining to this line of research.

\medskip

In order to keep the length of this paper within bounds, in Section 4 we will present applications that are related to a very specific setting, namely edge-counting in random geometric graphs with possibly large connections. The power and flexibility of the results proved in the present work are further illustrated in the companion paper \cite{LRP2}, where the following applications are developed in full detail: 

\begin{itemize}

\item[(i)] analytic bounds for the normal approximation of $U$-statistics based on marked point processes, in particular $U$-statistics with rescaled kernels;

\item[(ii)] bounds for general subgraph counting in the disk graph model under any regime;

\item[(iii)] an exhaustive characterization of the asymptotic behavior of geometric $U$-statistics; 

\item[(iv)] applications to the boolean model, and to subgraph counting in disk graph models with random radius.

\end{itemize}

The rest of this section is devoted to the formal presentation of the main problems that are addressed in this paper.

\subsection{Poisson measures}\label{ss:intropoisson}

Throughout the paper $(Z,\mathscr{Z},\mu) $ is a measure space such that
$Z$ is a Borel space, $\mathscr{Z}$ is the associated Borel $\sigma$-field, and $\mu$ is a $\sigma$-finite
Borel measure with no atoms. We write $\mathscr{Z}_{\mu} = \{ B\in \mathcal{Z}: \mu(B)< \infty \}$ to denote the subclass of $\mathscr{Z}$ composed of sets with finite measure. Also, we shall write
$\eta = \{\eta(B) : B\in \mathscr{Z}_{\mu} \} $ to indicate a {\it Poisson measure} on $(Z,\mathcal{Z}) $ with {\sl control} $\mu$. In other words, $\eta $ is a collection of random variables defined on some probability space $(\Omega, \mathscr{F}, P) $, indexed by
the elements of $\mathscr{Z}_{\mu} $ and such that: (i) for every $B,C \in \mathcal{Z}_{\mu}$ such that $B \cap C = \varnothing$, the random variables $ \eta(B)$ and $ \eta(C)$ are independent;  (ii) for every $B \in \mathscr{Z}_{\mu} $, $\eta(B)$ has a Poisson distribution with mean $\mu(B)$. We shall also write $\hat{\eta}(B) = \eta(B) - \mu(B)$, $B\in \mathscr{Z}_\mu$, and $\hat{\eta} = \{\hat{\eta}(B) : B\in \mathscr{Z}_{\mu} \}$. A random measure verifying property (i) is usually called ``completely random'' or ``independently scattered'' (see e.g. \cite{PeTa} for a general introduction to these concepts).  

\begin{rem}{\rm
As it is customary, by a slight abuse of notation, we shall often write $ x \in \eta$ in order to indicate that the point $x\in Z$ is charged by the random measure $\eta(\cdot)$.
}
\end{rem}

In this paper, we shall focus on sequences of random variables $\{F_n : n\geq 1\}$ having a finite Wiener-It\^o chaotic decomposition, that is, such that
\begin{equation}\label{e:fn}
F_n = \sum_{i=1}^k I_{q_i}(f_i^{(n)}), \quad n\geq 1,
\end{equation}
where the symbol $I_{q_i}$ indicates a multiple Wiener-I\^o integral of order $q_i$ with respect to $\hat{\eta}$, the integer $k$ does not depend on $n$, and each $f_i^{(n)}$ is a non-zero symmetric kernel from $Z^{q_i}$ to $\mathbb{R}$ (see Section \ref{ss:mwii} below for details). We will be specifically concerned with the forthcoming Problem \ref{p:chaos}. Recall that, given random variables $U,Y\in L^1(P)$, the {\it Wasserstein distance} between the law of $U$ and the law of $Y$ is defined as the quantity
\[
d_W(U,Y) = \sup_{f\in {\rm Lip}(1)} \big|E[f(U)] -E[f(Y)]\big|, 
\]
where ${\rm Lip}(1)$ indicates the class of Lipschitz real-valued function with Lipschitz constant $\leq 1$. It is well-known that the topology induced by $d_W$, on the class of probability measures on the real line, is strictly stronger than the one induced by the convergence in distribution.

\begin{prob}\label{p:chaos}{\rm Find analytic conditions on the kernels $\{f_i^{(n)}\}$ ensuring that the sequence
\[
\tilde{F}_n := \frac{F_n}{\sqrt{{\rm Var}(F_n)}}, \quad n\geq 1,
\]
converges in distribution, as $n\to \infty$, to a standard Gaussian random variable $N\sim \mathscr{N}(0,1)$, in the sense of the Wasserstein distance. Determine under which assumptions these conditions are also necessary, and find explicit upper bounds for the sequence
\[
d_W(\tilde{F}_n, N), \quad n\geq 1.
\]
}
\end{prob}

\medskip

We will deal with Problem \ref{p:chaos} in Section \ref{s:normalapp}, where it is shown that a convenient solution can be deduced by using {\it contraction operators}. Among other features, these operators provide a neat way to deal with the product of multiple stochastic integral, and virtually replace the use of diagram formulae -- see e.g. \cite{PeTa}. As anticipated, we will see that, in the specific case of random variables as in (\ref{e:fn}) such that $f_i^{(n)}\geq 0$, our results lead to necessary and sufficient conditions that are analogous to the so-called `fourth moment theorems' for sequences of multiple integrals in a Gaussian setting -- see \cite{nuapec}. 

\begin{rem}{\rm Problem \ref{p:chaos} is also explicitly studied in \cite[Section 4]{lesmathias}. In particular, Theorem 4.1 in \cite{lesmathias} provides bounds in the Wasserstein distance  for random variables having a finite chaotic decomposition, where the bounds are expressed in terms of expectations of inner products of multiple integral stochastic processes. On the other hand, Theorem 4.7 in \cite{lesmathias} provides an analytic bound, involving sums over partitions, for the normal approximation of absolutely convergent $U$-statistics. Here, we call `analytic bound' any upper bound only involving deterministic transformations of the kernel determining the $U$-statistic, without any additional probabilistic component.
}
\end{rem}

\subsection{Random graphs}

As anticipated, we shall now apply our main theoretical results to the study of geometric random graphs whose edge-counting statistics satisfy a CLT. The class of geometric random graphs considered below allow for long connections, in the sense that the geometric rule used to define edges is based on the use of arbitrarily large sets and therefore is not local. It is worth noting by now that our setting represents a natural generalization of the so called {\it Gilbert graphs} -- see Example \ref{ex:gilbert} below. Also, as explained in Remark \ref{rmk:stab} below, part of the models we consider cannot be dealt with by directly using the powerful theory of {\it stabilization} (see e.g. \cite{MolKen11}).

\medskip

Now let the notation introduced in the previous section prevail. In what follows, we shall denote by $W$ (as in `window') a measurable subset of $Z$ such that $\mu(W) <\infty$. We first introduce the notion of a {\it geometric random graph} based on the restriction of the Poisson measure $\eta$ to $W$, and on some symmetric set $H\subset Z\times Z$.

\begin{defi}[Random geometric graphs]\label{d:randomg}{\rm Let $H\subset Z\times Z$ be such that $\mu^2(H)<\infty$, $H$ is symmetric (that is, for every $(x,y) \in H$, one also has $(y,x)\in H$) and $H$ is non-diagonal  (that is, $H$ does not contain any pair of the type $(x,x)$). 
\begin{itemize}
\item[\rm (a)] The {\it random geometric graph} based on $\eta$, $W$ and $H$ is the undirected random graph \[
G =G(\eta,W,H) = (V,\mathcal{E}),
\] 
such that: (i) the vertices of $G$ are given by the class $V = \eta \cap W = \{x\in \eta : x\in W\}$, and (ii) a pair $\{x,y\}$ belongs to the set $\mathcal{E}$ of the edges of $G$ if and only if $(x,y)\in H$. We observe that, since $H$ is non-diagonal, $G$ has no loops, that is: $G$ does not contain any edge of the type $\{x,x\}$.

\item[\rm (b)] Assume in addition that $Z$ is a vector space. The random geometric graph at Point (a) is said to be {\it stationary} if there exists a set $\fact{H} \subset Z$ such that 
\[
H =\big \{(x,y)\in Z^2 : x-y\in \fact{H}\big \}.
\]
Note that, since $H$ is symmetric, one has necessarily that $\fact{H} = -\fact{H}$; moreover, since $H$ has no diagonal components, $0\notin \fact{H}$.
\end{itemize}
}
\end{defi}

\begin{example}\label{ex:gilbert}{\rm 
\begin{enumerate}
\item[\rm (i)] The class of random geometric graphs introduced above generalizes the notion of a {\it Gilbert graph}, obtained by taking $Z$ equal to some metric space (endowed with a distance $d$) and $H = \big\{(x_1,x_2)\in Z^2 : 0< d(x_1,x_2)< \delta\big\}$, $\delta >0$. In this case, the random geometric graph $G$ is obtained by connecting two vertices $v_1,v_2\in \eta \cap W$ if and only if $d(v_1,v_2) < \delta$. See e.g. \cite{penrosebook}. 

\item[\rm (ii)] If $Z = \mathbb{R}^d$ and $H = \big\{(x_1,x_2)\in Z^2 : 0< \|x_1 - x_2\|_{\mathbb{R}^d} < \delta\big\}$, then the corresponding geometric random graph is stationary with $\fact{H} = B(0,\delta)\backslash \{0\}$, where $B(0,\delta)\subset Z$ stands for the open ball of radius $\delta$ centered at the origin.  
Graphs of this type are customarily called {\it interval graphs} when $Z = \mathbb{R}$, and {\it disk graphs} when $Z=\mathbb{R}^2$ -- see e.g. \cite{CCJ, DFR, FV, McDiarmid, Muller}  for recent developments on the subject.
\end{enumerate}

}
\end{example}

In Section \ref{s:graphs}, we shall use our general results in order to deal with the following problem.

\begin{prob}\label{p:graphs}{\rm Fix a set $W\subset Z$, as well as a symmetric function $g : W\times W \to \mathbb{R}$. %and a measure $\alpha$ on $(Z^2,\mathscr{Z}^2)$ such that $\alpha(W\times W)<\infty$ 
Consider the following objects: (i) a collection of Poisson measures $\eta_\lambda$, $\lambda >0$, with $\sigma$-finite and non atomic control measures $\mu_\lambda$ such that $\mu_\lambda(W)<\infty$ and $\mu_\lambda(W)\uparrow \infty$, as $\lambda\to \infty$, and (ii) a collection of non-diagonal non-negligible symmetric sets $H_\lambda\subset Z\times Z$, $\lambda >0$. Characterize those classes $\{\mu_\lambda, \, H_\lambda : \lambda>0\}$ such that the random variables 
\begin{equation}\label{e:graphstat}
F(g,W; \eta_\lambda,H_\lambda) = \sum_{\{x,y\}\in \mathcal{E}_\lambda} g(x,y), \quad \lambda>0,
\end{equation}
(where $G_\lambda = (V_\lambda,\mathcal{E}_\lambda)$ is the geometric random graph based on $\eta_\lambda$, $W$ and $H_\lambda$) verify the limit relation
\begin{equation}\label{e:graphclt}
\widetilde{F}(g,W; \eta_\lambda,H_\lambda) :=\frac{F(g,W; \eta_\lambda,H_\lambda) - E[F(g,W; \eta_\lambda,H_\lambda)]}{\sqrt{{\rm Var}(F(g,W; \eta_\lambda,H_\lambda))}} \stackrel{\rm Law}{\longrightarrow}N\sim \mathscr{N}(0,1),
\end{equation}
as $\lambda\to \infty$.
}
\end{prob}

\begin{rem}{\rm

\begin{enumerate}

\item[(i)] At this stage, the role of the window $W$ might seem immaterial, and indeed the substance of Problem \ref{p:graphs} does not change if one takes $W=Z$. However, the above formulation allows for the more general case of a window $W = W_\lambda$ possibly depending on $\lambda$. Moving windows of this type appear in Section 4, as well as in the paper \cite{LRP2}.

\item[\rm (ii)] In many examples and applications, one considers sets $H_\lambda$ such that $\alpha(H_\lambda\cap (W\times W) )\downarrow 0$, as $\lambda \to \infty$, for some fixed measure $\alpha$ on $W\times W$. Heuristically, the fact that $\mu_\lambda(W) \uparrow \infty$  and $\alpha(H_\lambda\cap (W\times W) )\downarrow 0$ ensures that the following phenomenon takes place: as $\lambda$ grows, more and more vertices and edges are added to the geometric graph, whereas old edges are deleted as a consequence of the asymptotic negligibility of $H_\lambda\cap (W\times W)$. Solving Problem \ref{p:graphs} in this framework is equivalent to characterizing all sequences of random geometric graphs such that the addition of vertices and the cancellation of edges compensate, thus generating asymptotic Gaussian fluctuations.

\end{enumerate}
}
\end{rem}

When specialized to the case of Gilbert graphs on $Z = \mathbb{R}^d$, Problem \ref{p:graphs} is tackled in the classic reference \cite[Chapter 3]{penrosebook} as a special case of general subgraph counting. A comparison with the results of \cite[Chapter 3]{penrosebook} is provided in Section \ref{ss:exex} below. A complete solution of Problem \ref{p:graphs} for general subgraph counting in Gilbert graphs, based on the techniques developed in this paper, is presented in \cite[Section 3]{LRP2}. See also \cite[Section 6.2]{lesmathias}.

\begin{rem}
\label{rmk:stab}{\rm 
Assume that, for every $x\in \eta$, there exists a random radius $R_{x}$  such that all the $y$ connected to $x$ in the random graph lie in the ball with center $x$ and radius $R_{x}$. Then, the variable $F = F(1,W; \eta_\lambda,H_\lambda)$ in (\ref{e:graphstat}) is \emph{stabilizing}, meaning that $F$ can be written in the form
\begin{equation*}
F=\sum_{x\in \eta}\xi(x,\eta),
\end{equation*} where $\xi$ is such that $\xi(x,\eta)$ is not modified by adding or removing a finite number of points to $\eta$ outside the ball with center $x$ and radius $R_{x}$ (see \cite{MolKen11} for more details on this topic). In our case, to fit the framework of formula (\ref{e:graphstat}) in the case $g=1$, $\xi(x,\eta)$ should be defined as 
\begin{equation*}
\xi(x,\eta)=\frac{1}{2}\#(\{y\neq x:\, \{x,y\}\in \mathcal{E}_\lambda\}),
\end{equation*}
where $\# A$ indicates the cardinality of $A$. The CLTs presented for instance in \cite{BarPen05, PenYuk08} cover well this case. Remark that in this particular framework of a deterministic connection rule, stabilization theory only allows for a bounded length, while we consider here models where points can have arbitrarily long connections.}
\end{rem}

\medskip

The rest of the paper is organized as follows. In Section \ref{s:prel}, we discuss several background results concerning Poisson measures, Wiener chaos and $U$-statistics. Section \ref{s:normalapp} contains our main abstract results concerning the normal approximation of random variables having a finite chaotic decomposition. Section \ref{s:graphs} focuses on random graphs and on several analytical characterizations of associated CLTs. An Appendix (see Section \ref{APPENDIW}) provides some basic definitions and results of Malliavin calculus.

\section{Preparation}\label{s:prel}

\subsection{Multiple integrals and chaos}\label{ss:mwii}

As before, $(Z,\mathscr{Z},\mu) $ is a non-atomic Borel measure space, and $\eta$ is a Poisson measure on $Z$ with control $\mu$.

\begin{rem} \label{rmk1}
\rm{
By virtue of the assumptions on the space $(Z,\mathscr{Z},\mu)$, and to simplify the discussion,
we will assume throughout the paper that $(\Omega,\mathscr{F},P)$ and $\eta$ are such that
\[ \Omega = \left\{ \omega = \sum_{j=1}^{n} \delta_{z_j},n\in   \mathbb{N} \cup \{\infty\},z_j\in Z    \right\}, \]
where $\delta_z$ denotes the Dirac mass at $z$, and $\eta$ is defined as the
\textit{canonical mapping}
\[(\omega,B) \mapsto \eta(B)(\omega) = \omega(B) ,\quad B\in \mathscr{Z}_{\mu},\quad \omega\in\Omega.\]  Also, the $\sigma$-field $\mathscr{F}$ will be always supposed to be the $P$-completion of the $\sigma$-field generated by $\eta$. }
\end{rem}

Throughout the paper, for $p\in [1,\infty)$, the symbol $L^p(\mu)$ is shorthand for $L^p(Z,\mathscr{Z},\mu)$.
  For an integer $q\geq 2$, we shall write $L^p(\mu^q) := L^p(Z^q, \mathscr{Z}^{\otimes q}, \mu^{q}) $, whereas  $L^p_s(\mu^q)$ stands for the subspace of $L^p(\mu^q)$ composed of functions that are $\mu^{q}$-almost everywhere symmetric. Also, we adopt the convention $L^p(\mu) = L_s^p(\mu) =L^p(\mu^1) =L_s^p(\mu^1) $ and use the following standard notation: for every $q\geq 1$ and every $f,g\in L^2(\mu^q)$,
  $$ \langle f,g \rangle_{L^2(\mu^q)} = \int_{Z^q} f(z_1,...,z_q)g(z_1,...,z_q)\mu^q (dz_1,...,dz_q), \quad \|f\|_{L^2(\mu^q)} = \langle f,f \rangle^{1/2}_{L^2(\mu^q)} . $$

For every $f\in L^2(\mu^q)$, we denote by $\widetilde{f}$ the canonical symmetrization of $f$, that is,
      \[\widetilde{f}(x_1,\ldots,x_{q})=\cfrac{1}{q!}
\sum_\sigma f (x_{\sigma(1)},\ldots,x_{\sigma(q)}), \]
where $\sigma $ runs over the $q! $ permutations of the
set $\{1,\ldots,q \}$. Note that $
\|\tilde{f}\|_{L^2(\mu^q)} \leq  \|f\|_{L^2(\mu^q)}$ (to see this, use for instance the triangular inequality) .

%For every $f\in L^2_s(\mu^n) $, $n\geq 1$, and every fixed $z\in Z$, we write $f(z,\cdot)$ to indicate the function defined on $Z^{n-1}$ given by $(z_1,\ldots,z_{n-1}) \mapsto f(z,z_1,\ldots,z_{n-1})$. Accordingly, $\widetilde{f(z,\cdot)}$ stands for the symmetrization of the function $f(z,\cdot)$ (in $(n-1)$ variables). Note that, if $n=1$, then $f(z,\cdot)=f(z)$ is a constant.

\begin{defi}{\rm
For every deterministic function $h\in L^2(\mu)$, we write
\[I_1(h)=\hat{\eta}(h) = \int_Z h(z) \hat{\eta}(dz) \] to indicate the {\it Wiener-It\^o
integral} of $h$ with respect to $\hat{\eta}$. For every $q\geq 2$ and every $f\in L_s^2(\mu^q)$, we denote by $I_q(f)$
the {\it multiple Wiener-It\^o integral}, of order $q$, of $f$ with respect to $\hat{\eta}$. We also set $I_q(f)=I_q(\tilde{f})$, for every $f\in L^2(\mu^q)$ (not necessarily symmetric), and $I_0(b)=b$ for every real constant $b$.
}
\end{defi}

The reader is referred for instance to the monograph \cite{PeTa}, by Peccati and Taqqu, for a complete discussion of multiple Wiener-It\^o integrals and their properties (including the forthcoming Proposition \ref{P : MWIone} and Proposition \ref{P: MWIchaos}).

\begin{prop}\label{P : MWIone}
The following equalities hold for every $q,m\geq 1$, every $f\in L_s^2(\mu^q)$  and every $g\in L_s^2(\mu^m)$:
\begin{enumerate}
  \item[\rm 1.] $E[I_q(f)]=0$,
  \item[\rm 2.] $E[I_q(f) I_m(g)]= q!\langle f,g  \rangle_{L^2(\mu^q)} \1_{(q=m)} $
  {\rm (isometric property).}
\end{enumerate}
\end{prop}

The Hilbert space composed of the random variables with the form $I_q(f)$, where $q\geq 1$ and $f\in L^2_s(\mu^q)$, is called the $q$th \emph{Wiener chaos} associated with the Poisson measure $\eta$. The following well-known {\it chaotic representation property} is an essential feature of Poisson random measures. Recall that $\mathscr{F}$ is assumed to be generated by $\eta $.

\begin{prop}
[Wiener-It\^o chaotic decomposition] \label{P: MWIchaos} Every random variable \[F\in L^2(\Omega, \mathscr{F},P):=L^2(P)\]
admits a (unique) chaotic decomposition of the type
\begin{equation} \label{chaos}
F= E[F] + \sum_{i = 1}^{\infty} I_i(f_i),
\end{equation}
where the series converges in $L^2(P)$ and, for each $i\geq 1$, the kernel $f_i$ is an element
of $L^2_s(\mu^i)$.
\end{prop}

\subsection{Star contractions and multiplication formulae}\label{ss:star}
We shall now introduce {\it contraction operators}, and succinctly discuss some of their properties. As anticipated in the Introduction, these objects are at the core of our main results.

\smallskip

The kernel $f \star_r^l g$ on $Z^{p+q-r-l}$, associated with functions $f\in L^2_s(\mu^p) $ and $g \in L^2_s(\mu^q) $, where $p,q \geq 1$, $r=1,\ldots, p\wedge q$ and $l=1,\ldots,r $, is defined as follows:
\begin{eqnarray}
& & f \star_r^l
g(\gamma_1,\ldots,\gamma_{r-l},t_1,,\ldots,t_{p-r},s_1,,\ldots,s_{q-r}) \label{contraction} \\
&=& \int_{Z^l} \mu^l(dz_1,...,dz_l)
f(z_1,,\ldots,z_l,\gamma_1,\ldots,\gamma_{r-l},t_1,,\ldots,t_{p-r}) \nonumber \\
& & \quad\quad\quad\quad\quad\quad\quad\quad\quad\quad\quad\quad \times g(z_1,,\ldots,z_l,\gamma_1,\ldots,\gamma_{r-l},s_1,,\ldots,s_{q-r}). \nonumber
\end{eqnarray}
Roughly speaking, the star operator `$\,\star_r^l\,$' reduces the number of variables in the tensor product of $f$ and $g$ from $p+q$ to $p+q-r-l$: this operation is realized by first identifying $r$ variables in $f$ and $g$, and then by integrating out $l$ among them. To deal with the case $l=0$ for $r=0,\ldots, p\wedge q$, we set
\begin{eqnarray*}
& &f \star_r^0
g(\gamma_1,\ldots,\gamma_{r},t_1,,\ldots,t_{p-r},s_1,,\ldots,s_{q-r}) \\
&=& f(\gamma_1,\ldots,\gamma_{r},t_1,,\ldots,t_{p-r})
g(\gamma_1,\ldots,\gamma_{r},s_1,,\ldots,s_{q-r}),
\end{eqnarray*}
 and
$$ f \star_0^0 g (t_1,,\ldots,t_{p},s_1,,\ldots,s_{q}) =f\otimes g (t_1,,\ldots,t_{p},s_1,,\ldots,s_{q})= f(t_1,,\ldots,t_{p})
g(s_1,,\ldots,s_{q}). $$
By using the Cauchy-Schwarz inequality, one sees immediately that $f \star_r^r g$ is square-integrable for any choice of
$r=0,\ldots, p\wedge q$ , and every $f\in L^2_s(\mu^p) $, $g \in L^2_s(\mu^q) $. \\

\begin{rem}\label{r:fub}{\rm
For every
$1\leq p\leq q$ and every
$r=1,...,p$,
\begin{equation}\label{useful}
\int_{Z^{p+q-r}} (f\star_r^0 g)^2 d\mu^{p+q-r} = \int_{Z^r}
(f\star_p^{p-r} f)(g\star_q^{q-r} g) d\mu^r,
\end{equation}
for every $f\in L_s^2(\mu^p)$ and every $g\in L_s^2(\mu^q)$}
\end{rem}

The next statement contains an important {\it product formula} for Poisson multiple integrals (see e.g. \cite{PeTa} for a proof).
\begin{prop}
[Product formula] Let $f\in L^2_s(\mu^p) $ and $g\in
L^2_s(\mu^q)$, $p,q\geq 1 $, and suppose moreover that $f \star_r^l g
\in L^2(\mu^{p+q-r-l})$ for every $r=1,\ldots,p\wedge q $ and $
l=1,\dots,r$ such that $l\neq r $. Then,
\begin{equation} \label{e:product}
I_p(f)I_q(g) = \sum_{r=0}^{p\wedge q} r!
\left(
\begin{array}{c}
  p\\
  r\\
\end{array}
\right)
 \left(
\begin{array}{c}
  q\\
  r\\
\end{array}
\right)
 \sum_{l=0}^r
 \left(
\begin{array}{c}
  r\\
  l\\
\end{array}
\right)  I_{p+q-r-l} \left(\widetilde{f\star_r^l g}\right),
\end{equation}
 with the tilde $\sim$ indicating a symmetrization, that
is,
\[ \widetilde{f\star_r^l g}(x_1,\ldots,x_{p+q-r-l})=\cfrac{1}{(p+q-r-l)!}
\sum_\sigma f\star_r^l g(x_{\sigma(1)},\ldots,x_{\sigma(p+q-r-l)}), \]
where $\sigma $ runs over all $(p+q-r-l)! $ permutations of the
set $\{1,\ldots,p+q-r-l \}$.

\end{prop}

\subsection{About the Malliavin formalism}

For the rest of the paper, we shall use definitions and results related to Malliavin-type operators defined on the space of functionals of the Poisson measure $\eta$. Our formalism coincides with the one introduced by Nualart and Vives in \cite{nuaviv}. In particular, we shall denote by $D$, $\delta$, $L$ and $L^{-1}$, respectively, the {\it Malliavin derivative}, the {\it divergence operator}, the {\it Ornstein-Uhlenbeck generator} and its {\it pseudo-inverse}. The domains of $D$, $\delta$ and $L$ are written ${\rm dom} D$, ${\rm dom} \delta$ and ${\rm dom} L$. The domain of $L^{-1}$ is given by the subclass of $L^2(P)$ composed of centered random variables. For the convenience of the reader we have collected some crucial definitions and results in the Appendix (see Section \ref{APPENDIW}). Here, we just recall that, since the underlying probability space $\Omega$ is assumed to be the collection of discrete measures described in Remark \ref{rmk1}, then one can meaningfully define the random variable $\omega\mapsto F_z (\omega) =F(\omega + \delta_z),\, \omega \in \Omega, $  for every given random variable $F$ and every $z\in Z$, where $\delta_z$ is the Dirac mass at $z$. One can therefore prove that the following neat representation of $D$ as a {\it difference operator} is in order.
\begin{lemme}\label{L : diff}
For each $F\in {\rm dom} D$,
$$ D_z F = F_z - F ,\,\,  \text{a.e.-} \mu(dz). $$
\end{lemme}
A complete proof of Lemma \ref{L : diff} can be found in \cite{nuaviv}. 

\subsection{$U$-statistics}
Following \cite[Section 3.1]{lesmathias}, we now introduce the concept of a $U$-statistic associated with the Poisson measure $\eta$.

\begin{defi}[$U$-statistics] {\rm Fix $k\geq 1$. A random variable $F$ is called a $U$-{\it statistic of order} $k$, based on the Poisson measure $\eta$, if there exists a kernel $f\in L^1_s(\mu^k)$ such that
\begin{equation}\label{e:ustat}
F = \sum_{(x_1,...,x_k)\in \eta^k_{\neq}} f(x_1,...,x_k), 
\end{equation}
where the symbol $\eta^k_{\neq}$ indicates the class of all $k$-dimensional vectors $(x_1,...,x_k)$ such that $x_i\in \eta$ and $x_i\neq x_j$ for every $1\leq i\neq j \leq k$. As made clear in {\rm \cite[Definition 3.1]{lesmathias}}, the possibly infinite sum appearing in (\ref{e:ustat}) must be regarded as the $L^1(P)$ limit of objects of the type $\sum_{(x_1,...,x_k)\in \eta^k_{\neq}\cap A_n} f(x_1,...,x_k)$, $n\geq 1$, where the sets $A_n \in Z^k$ are such that $\mu^k(A_n)<\infty$ and $A_n \uparrow Z^k$, as $n\to \infty$.}
\end{defi}

Plainly, a $U$-statistic of order one is just a linear functional of $\eta$, with the form \[\sum_{x\in \eta}f(x)  = \int_Z f(x)\eta(dx),\] for some $f\in L^1(\mu)$. The following statement, based on the results proved  by Reitzner and Schulte in \cite{lesmathias}, collects two crucial properties of $U$-statistics.

\begin{thm}[See \cite{lesmathias}]\label{t:ustat} Let $F\in L^1(P)$ be a $U$-statistic as in {\rm(\ref{e:ustat})}. Then, the following two properties hold.
\begin{enumerate}

\item[\rm (a)] The expectation of $F$ is given by 
\begin{equation}\label{e:campbell}
E[F] = \int_{Z^k} f(z_1,...,z_k) \mu^k(dz_1,...,dz_k).
\end{equation}

\item[\rm (b)] If $F$ is also square-integrable, then necessarily $f\in L^2_s(\mu^k)$, and the Wiener-It\^o representation {\rm (\ref{chaos})} of $F$ is such that $f_i = 0$, for $i\geq k+1$, and
\begin{equation}\label{e:lastpenrose}
f_i(x_1,..,x_i) = \binom{k}{i}\int_{Z^{k-i}} f(x_1,...,x_i,z_{1},...,z_{k-i})\, \mu^{k-i}(dz_1,...,dz_{k-i})
\end{equation}
for every $i=1,...,k$. In particular, $f_k = f$. For each $i=1,\dots,k$, one has that $f_{i}\in L_{s}(\mu^i)\cap L_{s}^2(\mu^i)$.
\end{enumerate}
\end{thm}

One should note that formula (\ref{e:lastpenrose}) follows from an application of the results proved by Last and Penrose in \cite{lastpenrose}. 

\subsection{$U$-statistics and random graphs}

In this paper, we will be interested in characterizing the Gaussian fluctuations of $U$-statistics having a specific {\it support}. In particular this allows one to deal with the set of `local $U$-statistics'' introduced by Reitzner and Schulte in \cite[Section 6]{lesmathias}. Recall that a set $H\in Z^k$ is called {\it symmetric} if the following implication holds: if $(x_1,...,x_k)\in H$, then $(x_{\sigma(1)},...,x_{\sigma(k)})\in H$ for every permutation $\sigma$ of $\{1,...,k\}$.   

\begin{defi}[Support of a $U$-statistic]\label{d:glu}{\rm Let $k\geq 2$, and let $H\subset Z^k$ be a measurable symmetric set. A $U$-statistic $F$ as in (\ref{e:ustat}) is said to have {\it support in $H$} if the function $f$ is such that 
\[
f(x_1,...,x_k) = 0, \quad \forall (x_1,...,x_k) \notin H.
\]
}
\end{defi}

%\query{What is not a GLU?}

\begin{example}[Local $U$-statistics]{\rm Let $Z$ be a metric space. Then, the class of {\it local $U$-statistics}, as defined in \cite[Section 6]{lesmathias}, coincides with the family of $U$-statistics having support in a set of the type $H = \big\{(x_1,...,x_k) : {\rm diam}(\{x_1,...,x_k\}) < \delta\big\}$ for some $\delta >0$. Here, the symbol ${\rm diam}(B)$ is shorthand for the diameter of $B$.
}
\end{example}

We shall now point out a well-known connection between $U$-statistics and hypergraphs. Recall that a {\it hypergraph} of order $k\geq 2$ is a pair $(V,\mathcal{E})$, where $V= (v_1,...,v_m)$ is a set of vertices, and $\mathcal{E} = (E_1,...,E_s)$ is a collection of (possibly non-disjoint) subsets of $V$ (called {\it edges}), such that each $E_i$ contains exactly $k$ elements; in particular a hypergraph of order $2$ is an undirected graph. 

\begin{rem}[$U$-statistics as graph statistics]{\rm 
\begin{enumerate}

\item[\rm (i)] Let $k\geq 2$, let $F$ be a  $U$-statistic as in (\ref{e:ustat}), and assume that $f = {\bf 1}_{W^{k}} \times {\bf 1}_H$, where $W\subset Z$ is some set (usually called a `window') such that $\mu(W)<\infty$. Then, the random variable $\frac{1}{k!} F$ counts the number of edges in the random hypergraph $(V,\mathcal{E})$, obtained as follows: $V = \eta \cap W $, and the class of edges $\mathcal{E}$ is composed of all subsets $\{x_1,...,x_k\} \subset V$ such that $(x_1,...,x_k)\in H$.

\item[\rm (ii)] If $k=2$, $Z$ is some metric space (endowed with a distance $d$) and $H = \big\{(x_1,x_2)\in Z^2 : d(x_1,x_2)< \delta\big\}$ and  $f = {\bf 1}_W \times {\bf 1}_H$, then the random variable $\frac12 F$ counts the number of edges in the undirected graph whose vertices $V$ are given by the points of $W$ charged by $\eta$ and such that two vertices $v_1,v_2$ are connected by an edge if and only if $0<d(v_1,v_2) < \delta$. These are the `Gilbert random graphs' discussed in Example \ref{ex:gilbert}(i).
\end{enumerate}
}
\end{rem}

To conclude, we present the notion of a {\it stationary} $U$-statistic. It will play an important role in Section \ref{s:graphs}.

\begin{defi}[Stationary $U$-statistics]\label{d:stationary}{\rm Fix $k\geq 2$, assume that $Z$ is a vector space, and let $F$ be a $U$-statistic of the type (\ref{e:ustat}), having support in a symmetric set $H$). We shall say that $F$ is {\it stationary} if there exists $\fact{H}\subset  Z^{k-1}$ such that
\begin{equation}\label{e:homu}
H = \{(x_1,...,x_k) : (x_2-x_1,x_3-x_1,...,x_k-x_1)\in \fact{H}\}. 
\end{equation}
}
\end{defi}

\begin{example}{\rm Consider the case $k=2$, $Z = \mathbb{R}^d$ and $H = \big\{(x_1,x_2)\in Z^2 : \|x_1 - x_2\|_{\mathbb{R}^d} < \delta\big\}$. Then, the corresponding $U$-statistic $F$ is stationary, with $\fact{H} = B(0,\delta)$, where $B(0,\delta)\subset Z$ stands for the open ball of radius $\delta$ centered at the origin. See Example \ref{ex:gilbert}(ii).
}
\end{example}

\section{Normal approximations for finite chaotic expansions}\label{s:normalapp}

\subsection{Framework}

We shall tackle Problem \ref{p:chaos}, by focussing on the normal approximation of random variables $F$ having the form
\begin{equation}\label{e:genrv}
F = E[F]  + \sum_{i=1}^k I_{q_i}(f_i),
\end{equation} 
where:
\begin{itemize}

\item[--] $k\geq 1$ is an integer;

\item[--] the integers $q_i$, $i=1,...,k$, are such that $1\leq q_1 <q_2 <\cdots <q_k$;

\item[--] the symbol $I_q$ indicates a multiple Wiener-It\^o integral of order $q$, with respect to a centered Poisson measure $\hat{\eta} = \eta -\mu$, where $\eta$ is a Poisson measure on the Borel measurable space $(Z,\mathscr{Z})$, with deterministic and $\sigma$-finite control measure $\mu$;

\item[--] each kernel $f_i$ is a nonzero element of $ L^2_s( \mu^{q_i})$, and the class $\{f_i : i=1,...,k\}$ verifies in addition the forthcoming Assumption \ref{a:tech}.

\end{itemize}

\begin{as}[Technical assumptions on integrands]\label{a:tech}{\rm Let the notation of Section \ref{ss:star} prevail. Every random variable of the type (\ref{e:genrv}) considered in the sequel of this paper is such that the following properties (i)-(iii) are verified.
\begin{enumerate}

\item[(i)] For every $i=1,...,d$ and every $r=1,...q_i$, the kernel $f_i\star_{q_i}^{q_i-r} f_i $ is an element of $L^2(\mu^{r})$.

\item[(ii)] For every $i$ such that $q_i\geq 2$, every contraction of the type $(z_1,...,z_{2q_i - r- l})\mapsto |f_i|\star_r^l |f_i| (z_1,...,z_{2q_i - r- l})$ is well-defined and finite for every $r=1,...,q_i$, every $l=1,...,r$ and every $(z_1,...,z_{2q_i - r- l})\in Z^{2q_i-r-l}$.

\item[(iii)] For every $ i,j=1,...,d$ such that $\max(q_i,q_j) >1$, for every $k = |q_i - q_j| \vee 1,..., q_i+q_j-2$ and every $(r,l)$ verifying $k = q_i+q_j -2-r-l$, 
\[
\int_Z \left[\sqrt{ \int_{Z^k} (f_i(z,\cdot)\star_r^l f_j(z,\cdot))^2 \,\,d\mu^k  }\,\,\,\right]\mu(dz)<\infty,
\]
where, for every fixed $z\in Z$, the symbol $f_i(z,\cdot)$ denotes the mapping $(z_1,...,z_{q-1}) \mapsto f_i(z,z_1,...,z_{q-1})$.
\end{enumerate}

}
\end{as}

\begin{rem}{\rm According to \cite[Lemma 2.9 and Remark 2.10]{PecZheng}, Point (i) in Assumption \ref{a:tech} implies that the following properties (a)-(c) are verified:

\begin{enumerate}

\item[(a)] for every $1\leq i<j\leq k$, for every $r=1,...,q_i\wedge q_j$ and every $l=1,...,r$, the contraction $f_i \star_r^l f_j$ is a well-defined element of $L^2(\mu^{q_i+q_j-r-l})$;

\item[(b)] for every $1\leq i\leq j\leq k$ and every $r=1,...,q_i$, $f_i\star_r^0 f_j$ is an element of $L^2(\mu^{q_i+q_j-r})$;

\item[(c)] for every $i=1,...,k$, for every $r=1,...,q_i$, and every $l=1,...,r\wedge (q_i-1)$, the kernel $f_i\star_r^l f_i $ is a well-defined element of $L^2(\mu^{2q_i-r-l})$.
\end{enumerate}

In particular, the multiplication formula (\ref{e:product}) implies that every random variable $F$ verifying Assumption \ref{a:tech} is such that $I_{q_i}(f_i)^2 \in L^2(P)$ for every $i=1,...,k$, yielding in turn that $E[F^4]<\infty$. Following a similar route, one can also show that, under Assumption \ref{a:tech}, the random variable $\langle DF, -DL^{-1}F\rangle_{L^2(\mu)}$ is square-integrable (and not merely an element of $L^1(P)$).
}
\end{rem}

\begin{rem}{\rm Assumption \ref{a:tech} imply that Assumptions A-B-C in \cite{PecZheng} are verified, so that the computations therein can be directly applied in our framework.
}
\end{rem}

\begin{rem}{\rm For instance, Assumption \ref{a:tech} is verified whenever each $f_i$ is a bounded function with support in a rectangle of the type $B\times\cdots\times B$, where $\mu(B)<\infty$. 
%Many examples of random variables verifying Assumption \ref{a:tech} are discussed in Sections \ref{s:app1}-\ref{s:app2}.
}
\end{rem}

\subsection{A general bound}

Let $F$ be a random variable as in (\ref{e:genrv}) such that $E[F^2] =\sigma^2>0$ ($\sigma>0$) and $E[F] =m\in \R$, and consider a Gaussian random variable $N\sim \mathscr{N}(m, \sigma^2)$ with the same mean and variance. Then, a slight modification of \cite[Theorem 3.1]{PSTU} (the modification resides in the fact that we consider an arbitrary variance $\sigma^2$) yields the following estimates:
\begin{equation}\label{e:pstubound}
d_W(F,N) \leq B_1(F;\sigma) \leq B_2(F;\sigma), 
\end{equation}
where
\begin{eqnarray}
 B_1(F;\sigma)&=& \frac1\sigma E\left[ |\sigma^2 - \langle DF, -DL^{-1}F\rangle_{L^2(\mu)}|\right]\label{e:B1} \\ &&\quad\quad\quad\quad\quad+ \frac{1}{\sigma^2} \int_Z E[(D_zF)^2|D_zL^{-1}F|] \mu(dz),\notag\\
  B_2(F;\sigma) &=& \frac1\sigma\sqrt{E\left[ (\sigma^2 - \langle DF, -DL^{-1}F\rangle_{L^2(\mu)})^2\right]}\label{e:B2} \\
  &&\quad \quad\quad\quad\quad+ \frac{1}{\sigma^2} \int_Z E[(D_zF)^2|D_zL^{-1}F|] \mu(dz).\notag
\end{eqnarray}

The next statement shows that $B_2(F;\sigma^2)$ can be further bounded in terms of the contractions introduced in Section \ref{ss:star}.

\begin{thm}\label{t:main1} Let $F$ and $N$ be the random variables appearing in {\rm (\ref{e:pstubound})-(\ref{e:B2})}. Then, there exists a universal constant $C_0 = C_0(q_1,...,q_k)\in (0,\infty)$, depending uniquely on $q_1,...,q_k$, such that
\begin{equation}\label{e:newbound}
B_2(F;\sigma) \leq C_0\times B_3(F;\sigma),
\end{equation}
where 
\begin{eqnarray}\label{e:B3}
&& B_3(F;\sigma) \\
&& \notag = \frac{1}{\sigma}\left\{ \max_1\|f_i\star_r^l f_i\|_{L^2(\mu^{2q_i-r-l})} + \max_2 \|f_i\star_r^l f_j\|_{L^2(\mu^{q_i+q_j-r-l})} + \max_{i=1,...,k} \|f_i\|^2_{L^4(\mu^{q_i})}\right\}.
\end{eqnarray}
In the previous expression, $\displaystyle{\max_1}$ ranges over all $1\leq i  \leq k$ such that $q_i>1$, and all pairs $(r,l)$ such that $r\in\{1,...,q_i\}$ and $1\leq l\leq r\wedge (q_i-1)$, whereas $\displaystyle{\max_2}$ ranges over all $1\leq i < j \leq d$ and all pairs $(r,l)$ such that $r\in\{1,...,q_i\}$ and $l\in\{1,...,r\}$. When $q_1 = 1$, one can replace $ \| f_{1} \|^2_{L^4(\mu)}$ in the previous bound by the smaller quantity
\begin{equation}\label{e:rep}
\frac{\| f_1\|^3_{L^3(\mu)}}{\| f_1\|_{L^2(\mu)}}.
\end{equation}
\end{thm}

\noindent{\it Proof of Theorem {\rm \ref{t:main1}}}. Without loss of generality, we can assume that $m=0$. Also, throughout this proof, we write $F_i =I_{q_i}(f_i)$ and $\sigma_i^2 = E[I_{q_i}(f_i)^2] = q_i ! \| f_i\|^2_{L^2(\mu^{q_i})}$, $i=1,...,k$, in such a way that $\sigma^2 = E[F^2] =\sum_{i=1}^kE[F_i^2] =\sum_{i=1}^k\sigma_i^2$.  Now write 
\begin{eqnarray*}
\sqrt{E\left[ (\sigma^2 - \langle DF, -DL^{-1}F\rangle_{L^2(\mu)})^2\right]}&\leq& \sum_{i=1}^k{\bf 1}_{q_i>1} \sqrt{E\left[ (\sigma_i^2 - \langle DF_i, -DL^{-1}F_i\rangle_{L^2(\mu)})^2\right]} \\
&+& \sum_{1\leq i\neq j\leq k} \sqrt{E\left[ (\langle DF_i, -DL^{-1}F_j\rangle_{L^2(\mu)})^2\right]},
\end{eqnarray*}
so that one can directly apply \cite[Proposition 5.5]{PecZheng} and deduce that there exists a constant $a=a(q_1,...,q_k)$ such that
\begin{eqnarray*}
&&\sqrt{E\left[ (\sigma^2 - \langle DF, -DL^{-1}F\rangle_{L^2(\mu)})^2\right]}\\&&\leq a\left\{ \max_1\|f_i\star_r^l f_i\|_{L^2(\mu^{2q_i-r-l})} + \max_2 \|f_i\star_r^l f_j\|_{L^2(\mu^{q_i+q_j-r-l})} \right\}.
\end{eqnarray*}
To conclude, observe that
\[
 \int_Z E[(D_zF)^2|D_zL^{-1}F|] \mu(dz) \leq  \int_Z E\left[\left(\sum_{i=1}^k | D_zF_i |\right)^2 \left(\sum_{i=1}^k|D_zL^{-1}F_i|\right)\right]\mu(dz), 
\]
so that \cite[Proposition 5.6]{PecZheng} implies that there exists a constant $b = b(q_1,...,q_k)$ such that
\[
 \int_Z E[(D_zF)^2|D_zL^{-1}F|] \mu(dz)\leq b\sigma \left\{ \max_1\|f_i\star_r^l f_i\|_{L^2(\mu^{2q_i-r-l})}  + \max_{i=1,...,k} \|f_i\|^2_{L^4(\mu^{q_i})}\right\}.
\]
Taking $C_0 =a+b$ yields the desired conclusion. The last assertion in the statement comes from the fact that, when $q_1 =1$, $\int_Z E[(D_zF)^2|D_zL^{-1}F|] \mu(dz) = \|f_1\|^3_{L^3(\mu)}$.
\qed

\begin{rem}{\rm According to \cite[Lemma 2.9]{PecZheng}, for every quadruple $(i,j,r,l)$ entering the expression of $\displaystyle{\max_2}$ in (\ref{e:B3}), the following estimate holds:
\begin{equation}\label{e:pzestimate}
\|f_i\star_r^l f_j\|^2_{L^2(\mu^{q_i+q_j-r-l})} \leq \|f_i\star_r^l f_i\|_{L^2(\mu^{2q_i-r-l})}\times \|f_j\star_r^l f_j\|_{L^2(\mu^{2q_j-r-l})} .
\end{equation}

}
\end{rem}

\medskip

\subsection{Estimates for positive kernels}\label{ss:d=1}

We shall now specialize Theorem \ref{t:main1} to the case on random variables having the form (\ref{e:genrv}) and such that $f_i \geq 0$. In particular, we shall unveil some useful connections between the quantity $B_3(F;\sigma)$ and the fourth cumulant of $F$.

\begin{rem}\label{r:pos}{\rm Random variables admitting a Wiener-It\^o chaotic expansion with positive kernels appear rather naturally in stochastic geometry. For instance, an application of (\ref{e:lastpenrose}) shows that any $U$-statistic with a positive kernel admits a Wiener-It\^o chaotic expansion of this type. Note that many {\it counting statistics} have the form of $U$-statistics with an integer-valued (and therefore nonnegative) kernel -- such as for instance the subgraph-counting statistics in random geometric graphs (see e.g. \cite[ Chapter 3]{penrosebook} and the references therein), or the statistics associated with hyperplane tessellations considered in \cite{hss}. 

}
\end{rem}

The following statement concerns random variables of the form (\ref{e:genrv}), in the special case where $E[F] =0$, $k=1$ and the multiple stochastic integral has a nonnegative kernel.

\begin{prop}[Fourth moment bound, I]\label{p:1bnd} Consider a random variable $F$ as in {\rm (\ref{e:genrv})}, with the special form $F= I_q(f)$, $q\geq 1$, where $f\geq 0$ and $E[F^2] = q!\|f\|^2_{L^2(\mu^q)} = \sigma^2>0$. Then, 
\[
B_3(F;\sigma) = \frac{1}{\sigma}\left\{ \max_1\|f\star_r^l f\|_{L^2(\mu^{2q-r-l})}+ \|f\|^2_{L^4(\mu^q)}\right\},
\]
and there exist universal constants $c_1 = c_1(q)< C_1 = C_1(q)$, depending uniquely on $q$, such that
\begin{equation}\label{e:bound1}
c_1\times B_3(F;\sigma)\leq \sqrt{E[F^4] - 3\sigma^4} \leq C_1\times B_3(F;\sigma).
\end{equation}
\end{prop}
\noindent{\it Proof}. Using the multiplication formula (\ref{e:product}), together with (\ref{useful}) and the positivity assumptions on $f$, one sees that $F^2$ can be written in the following form
\[
F^2 = \sigma^2 +I_{2q}(f\star_0^0 f) +R,
\]
where $R$ is a random variable orthogonal to the constants and to $I_{2q}(f\star_0^0 f)$ and such that
\[
  J \left\{ \max_1\|f\star_r^l f\|^2_{L^2(\mu^{2q-r-l})}+ \|f\|^4_{L^4(\mu^q)}\right\} \! \leq E[R^2]\! \leq\! K \left\{ \max_1\|f\star_r^l f\|^2_{L^2(\mu^{2q-r-l})}+ \|f\|^4_{L^4(\mu^q)}\right\},
\]
for some universal positive constants $J,K$ depending uniquely on $q$. The conclusion is obtained by using the relation
\[
E[I_{2q}(f\star_0^0 f)^2] = (2q)!\| \widetilde{f\star_0^0 f}\|^2_{L^2(\mu^{2q})} = 2\sigma^4 +\sum_{p=1}^{q-1}\frac{q!^4}{p!^2(q-p)!^2}  !\| \widetilde{f\star_p^p f}\|^2_{L^2(\mu^{2q-2p})},
\]
where we have used \cite[formula (11.6.30)]{PeTa}.
\qed

\medskip

The following general bound deals with random variables of the form (\ref{e:genrv}) and with positive kernels.

\begin{prop} [Fourth moment bound, II]\label{p:2bnd} Let $F$ be as in {\rm (\ref{e:genrv})}, with $k\geq 1$, and assume moreover that $E[F]=0$, and $E[F^2]=\sigma^2>0$, and $f_i\geq 0$ for every $i$. Then, there exists a universal constant $C_2 = C_2(q_1,...,q_k)$, depending uniquely on $q_1,...,q_k$, such that
\begin{equation}\label{e:4mb}
\left\{ \max_1\|f_i \star_r^l f_i \|_{L^2(\mu^{2q_i-r-l})}+ \max_{i=1,...,k}\|f_i\|^2_{L^4(\mu^{q_i})}\right\} \leq C_2 \sqrt{E[F^4] - 3\sigma^4}.
\end{equation}
\end{prop}
{\it Proof}. Write as before $F_i = I_{q_i}(f_i)$, and $\sigma_i^2 = E[F_i^2]$, $i=1,...,k$. We can now write
\begin{eqnarray*}
E[F^4] - 3\sigma^4 &=& \sum_{i=1}^k \{E[F_i^4] - 3\sigma_i^4\}\\
& +& 6\sum_{1\leq i<j\leq k} \{E[F_i^2 F_j^2] - \sigma_i^2\sigma_j^2\} +\sum_{(i_1,i_2,i_3,i_4)\in V_d} E[F_{i_1}F_{i_2}F_{i_3}F_{i_4}]\\
&:=& W+Y+Z,
\end{eqnarray*}
where $V_d$ stands for the collection of those $(i_1, i_2, i_3, i_4) \in  \{1, ..., k\}^4$, such that one of the following conditions is satisfied:
(a) $i_1 \neq  i_2 = i_3 = i_4$, (b) $i_1 \neq  i_2 = i_3\neq  i_4$ and $i_4\neq i_1$, (c) the elements of $(i_1,i_2,i_3, i_4)$ are all distinct. Applying the multiplication formula (\ref{e:product}) and exploiting the fact that each $f_i$ is nonnegative, we immediately deduce that $Y\geq 0$ and $Z\geq 0$, so that the desired conclusion follows from Proposition \ref{p:1bnd}.
\qed

\subsection{Conditions for asymptotic Gaussianity}\label{ss:nsclt}

This section contains a general statement (Theorem \ref{t:main2}) about the normal approximation of random variables admitting a finite chaotic decomposition. The first part of such a result provides sufficient conditions for Central Limit Theorems, that are directly based on Theorem \ref{t:main1}. As indicated in the subsequent parts, these conditions become necessary whenever the involved kernels are nonnegative. Theorem \ref{t:main2} is one of the main results of the paper, and is the main tool used to deduce the CLTs appearing in Section 4 and in \cite{LRP2}.

More precisely, in what follows we shall fix integers $k\geq 1$ and $1\leq q_1<q_2<...<q_k$ (not depending on $n$), and consider a sequence $\{F^{(n)} : n\geq 1\}$ of random variables with the form
\begin{equation}\label{e:genseq}
F^{(n)} =\sum_{i=1}^k I_{q_i}(f_i^{(n)}), \quad n\geq 1
\end{equation}
each verifying the same assumptions as the random variable $F$ appearing in {\rm (\ref{e:genrv})} (in particular, Assumption \ref{a:tech} is satisfied for each $n$). We also use the following additional notation: (i) $\sigma^2(n) = E[(F^{(n)})^2]$, (ii) $F^{(n)}_i =   I_{q_i}(f_i^{(n)})$, $i=1,...,d$, and (iii) $\sigma_i^2(n) = E[(F_i^{(n)})^2]$.  

\begin{rem}{\rm
One of the main achievements of the forthcoming Theorem \ref{t:main2} is the `fourth moment theorem' appearing at Point 3 in the statement, which only holds for random variables such that the kernels in the chaotic decomposition are nonnegative. As first proved in \cite{nuapec} (see also \cite[Chapter 5]{np-book}) an analogous result holds for general sequences of multiple Wiener-It\^o integrals with respect to a Gaussian process. In particular, in the Gaussian framework one does not need to assume that the integrands have a constant sign. Proving analogous statements in a Poisson setting is quite a demanding task, one of main reasons being the rather intricate multiplication formula (\ref{e:product}). Some previous partial findings in the Poisson case can be found in Peccati and Taqqu \cite{PeTaBer} (for sequences of double integrals) and in Peccati and Zheng \cite{PecZheng2011} (for sequences of multiple integrals having the form of homogeneous sums).
}
\end{rem}

\begin{rem}{\rm In the statement of Theorem \ref{t:main2}, we implicitly allow that the underlying Poisson measure $\eta$ also changes with $n$. In particular, one can assume that the associated control measure $\mu = \mu_n$ explicitly depends on $n$. This general framework is needed for the geometric applications developed in Section \ref{s:graphs}.
}
\end{rem}

\begin{thm}\label{t:main2} Let $\{F^{(n)} \}$ be a sequence of random variables as in {\rm (\ref{e:genseq})}, and assume that there exists $\sigma^2>0$ such that $\lim_{n\to \infty}\sigma^2(n) = \sigma^2$. Let $N\sim\mathscr{N}(0, \sigma^2)$.

\begin{enumerate}

\item[\rm 1.] For every $n$, one has the estimate
\begin{equation}\label{e:bnd}
d_W(F^{(n)},N)\leq C_0\times B_3(F^{(n)};\sigma(n)) +\frac{\sqrt{2/\pi}}{\sigma(n)\vee \sigma} |\sigma^2(n) - \sigma^2|.
\end{equation} 
In particular, if $B_3(F^{(n)};\sigma(n))\to 0$, as $n\to \infty$, then $d_W(F^{(n)},N)\to 0$ and therefore $F^{(n)} \stackrel{\rm Law}{\to} N$.

\item[\rm 2.] Assume that $f_i^{(n)}\geq 0$ for every $i,n$. Then, a sufficient condition in order to have that $B_3(F^{(n)};\sigma(n))\to 0$ is that $E[(F^{(n)})^4] - 3\sigma^4(n)\to 0$.

\item[\rm 3.] Assume that $f_i^{(n)}\geq 0$ for every $i,n$, and also that the sequence $(F^{(n)})^4$, $n\geq 1$, is uniformly integrable. Then, the following conditions {\rm (a)--(c)} are equivalent, as $n\to \infty$: {\rm (a)} $d_W(F^{(n)},N)\to 0$, {\rm (b)}  $B_3(F^{(n)};\sigma(n))\to 0$, and {\rm (c) }$E[(F^{(n)})^4] - 3\sigma^4(n)\to 0$.

\end{enumerate}

\end{thm}
{\it Proof}. 1. Let $N_n\sim\mathscr{N}(0, \sigma^2(n))$, $n\geq 1$. Then, one has that (see e.g. \cite[Proposition 3.6.1]{np-book})
\[
d_W(N_n,N) \leq \frac{\sqrt{2/\pi}}{\sigma(n)\vee \sigma} |\sigma^2(n) - \sigma^2|,
\]
so that the desired conclusion follows from Theorem \ref{t:main1} as well as the inequality \[d_W(F^{(n)},N)\leq d_W(F^{(n)},N_n)+d_W(N_n,N).\]

\noindent 2. If $E[(F^{(n)})^4] - 3\sigma^4(n)\to 0$, then relation (\ref{e:4mb}) implies that 
\[
\left\{ \max_1\|f^{(n)}_i \star_r^l f^{(n)}_i \|_{L^2(\mu^{2q_i-r-l})}+ \max_{i=1,...,k}\|f^{(n)}_i\|^2_{L^4(\mu^{q_i})}\right\}\to 0.
\]
Using (\ref{e:pzestimate}), we see that the last relation implies that $B_3(F^{(n)};\sigma(n))\to 0$, so that the desired conclusion follows from Point 1 in the statement.

\noindent 3. In view of Point 1 and Point 2 in the statement, we shall only prove that (a) $\Rightarrow$ (c). To prove this implication, just observe that $F^{(n)} \stackrel{\rm Law}{\to} N$ and $\{(F^{(n)})^4\}$ is uniformly integrable, then necessarily $E[(F^{(n)})^4] \to E[N^4] =  3\sigma^4$, so that the conclusion follows from the fact that $\sigma^2(n) \to \sigma^2$.
\qed

\begin{rem}{\rm
 A sufficient condition in order to have that the sequence $\{(F^{(n)})^4\}$ is uniformly integrable is the following: there exists some $\epsilon >0$ such that 
\[
\sup_{n\geq 1} E|F^{(n)}|^{4+\epsilon} <\infty.
\]
}
\end{rem}

\medskip

We shall use some estimates taken from \cite{PecZheng} (see, in particular, Table 2, p. 1505, therein). Given a thrice differentiable function $\phi : \R^k \to \R$, we set
$$  \|\phi''\|_{\infty}=\max\limits_{1\leq i_1 \leq i_2 \leq k} \sup\limits_{x\in \R ^k} \left| \cfrac{\partial^2}{\partial x_{i_1}  \partial x_{i_2}} \phi(x)\right| .$$
$$  \|\phi'''\|_{\infty}=\max\limits_{1\leq i_1 \leq i_2  \leq i_3\leq k} \sup\limits_{x\in \R^k} \left| \cfrac{\partial^3}{\partial x_{i_1} \partial x_{i_2} \partial x_{i_3}} \phi(x)\right| .$$

\begin{prop}\label{p:sumtosingle} Let the assumptions and notation of Theorem \ref{t:main2} prevail, and suppose that $B_3(F^{(n)};\sigma)\to 0$, as $n\rightarrow\infty$. Let $N_{n,i}\sim \mathscr{N}(0, \sigma^2_i(n))$, $i=1,...,k,$ Then, for every thrice differentiable function $\phi : \R^k \to \R$, such that $\|\phi''\|_{\infty},\|\phi'''\|_{\infty}<\infty$, one has that
\[
E[ \phi(F_1^{(n)},...,F_k^{(n)})] - E[ \phi(N_{n,1},...,N_{n,k})] \longrightarrow 0, \quad n\to \infty.
\]
\end{prop}
{\it Proof}. According to \cite{PecZheng}, the following estimate takes place:
\begin{eqnarray*}
&& | E[ \phi(F_1^{(n)},...,F_k^{(n)})] - E[ \phi(N_{n,1},...,N_{n,k})] |\\
&& \leq \frac{k}{2} \|\phi''\|_\infty \sqrt{\sum_{i,j=1}^k E[(\sigma^2_i(n){\bf 1}_{i=j} -\langle DF_i, -DL^{-1}F_j\rangle_{L^2(\mu)})^2] }\\
&+&\frac14 \|\phi'''\|_\infty \int_Z E\left[\left(\sum_{i=1}^k | D_zF_i |\right)^2 \left(\sum_{i=1}^k|D_zL^{-1}F_i|\right)\right]\mu(dz),
\end{eqnarray*}
so that the conclusion follows from (\ref{e:newbound}).

\qed

%\underline{ \bf J'ai redig\'e  jusqu'\` l\`a}

\section{Edge-counting in random geometric graphs: from Gaussian fluctuations to clustering}\label{s:graphs}

\subsection{Framework}\label{ss:frame}

Our aim is now to tackle Problem \ref{p:graphs}. Throughout this section we shall work under the following slightly more restrictive setting (we use the notation of Problem \ref{p:graphs}).

\begin{enumerate}
\item[--] For every $\lambda >0$, the control measure of $\eta_\lambda$ is $\mu_\lambda(dx) = \lambda\times \theta(dx)$, where $\theta$ is a $\sigma$-finite non-atomic measure on $(Z,\mathscr{Z})$.

%\item[--] $\alpha(dx,dy) = \theta(dx)\theta(dy) = \theta^2(dx,dy)$.

\item[--] The symmetric function $g : W\times W\to \mathbb{R}$ is bounded (this assumption can be relaxed -- see the discussion below). 

\end{enumerate}

In the forthcoming Section \ref{ss:genbounds}, we will show that the normal approximation of the random variables $F(g,W; \eta_\lambda,H_\lambda)$ (as defined in (\ref{e:graphstat})) can be completely characterized in terms of some diagonal restrictions of Cartesian products of the sets $H_\lambda$. Among several consequences, this remarkable phenomenon implicitly provides a new geometric interpretation of contraction operators.

\begin{defi}[Sets with diagonal restrictions]
\label{def:ABCDE}
 {\rm Given the sets $H_\lambda$, $\lambda>0$, defining the random variables $F(g,W; \eta_\lambda,H_\lambda)$ and $\widetilde{F}(g,W; \eta_\lambda,H_\lambda)$, respectively in (\ref{e:graphstat}) and (\ref{e:graphclt}), we define the following four sets

\begin{enumerate}
\item[--] $H^{(1)}_\lambda := \{(x_1,x_2,x_3)\in Z^3 : (x_1,x_2)\in H_{\lambda}, \,\, (x_1,x_3)\in H_{\lambda}\}\subset Z^3$;

\item[--] $H^{(2)}_\lambda := \{(x_1, x_2,x_3,x_4,x_5)\in Z^5 : (x_1,x_i)\in H_{\lambda}, \,\, \forall i=2,...,5\}\subset Z^5$;

%\item[--] $H^{(3)}_\lambda = \{(x_1,x_2,x_3,x_4)\in Z^4 : (x_1,x_2)\in H_\lambda, (x_3,x_4)\in H_\lambda, (x_1,x_4)\in H_\lambda, (x_3,x_2)\in H_\lambda\}\subset Z^4$;
\item[--] $H^{(3)}_\lambda = \{(x_1,x_2,x_3,x_4)\in Z^4 : (x_1,x_2)\in H_\lambda, (x_2,x_3)\in H_\lambda, (x_3,x_4)\in H_\lambda, (x_4,x_1)\in H_\lambda\}\subset Z^4$;

%\item[--] $H^{(4)}_\lambda := \{(x_1, x_2,x_3,x_4,x_5)\in Z^5 : (x_1,x_2)\in H_\lambda, (x_1,x_3)\in H_\lambda, (x_2,x_4)\in H_\lambda, (x_3,x_5)\in H_\lambda\}\subset Z^5$.
\item[--] $H^{(4)}_\lambda := \{(x_1, x_2,x_3,x_4,x_5)\in Z^5 : (x_1,x_2)\in H_\lambda, (x_2,x_3)\in H_\lambda, (x_3,x_4)\in H_\lambda, (x_4,x_5)\in H_\lambda\}\subset Z^5$.

\end{enumerate}
}
\end{defi}

\medskip

We shall also use the following quantities, defined for every $\lambda >0$:

\begin{enumerate}

\item[\rm (i) ] 
\begin{eqnarray*}
&& V^2_\lambda(g) =V^2_{1,\lambda}(g) +V^2_{2,\lambda}(g)\\ &&:= 4  \lambda^3 \int_{W^3\cap H^{(1)}_\lambda} g\star_1^0g(x_1,x_2,x_3)\theta^3(dx_1, dx_2,dx_3) + 2 \lambda^2\int_{W^2\cap H_\lambda} g^2(x_1,x_2)\theta^2(dx_1,dx_2);
\end{eqnarray*}

\item[\rm (ii) ]  $A_{\lambda}(g) =  \lambda^{5/2} \sqrt{\int_{W^5\cap H_\lambda^{(2)} }\prod_{i=2}^5 g(x_1,x_i) \theta^5(dx_1,...,dx_5)}$;

\item[\rm (iii) ] $B_{\lambda}(g) =  \lambda^{2} \sqrt{\int_{W^4\cap H_\lambda^{(3)} }g(x_1,x_2)g(x_3,x_4)g(x_1,x_4)g(x_2,x_3) \theta^4(dx_1,...,dx_4)}$;

\item[\rm (iv) ] $C_{\lambda}(g) =  \lambda^{3/2} \sqrt{\int_{W^3\cap H_\lambda^{(1)} }g^2(x_1,x_2)g^2(x_1,x_3) \theta^3(dx_1,dx_2,dx_3)}$;

\item[\rm (v) ] $D_{\lambda}(g) =  \lambda\sqrt{\int_{W^2\cap H_\lambda }g^4(x_1,x_2) \theta^2(dx_1,dx_2)}$;

\item[\rm (vi)] $E_{\lambda}(g) =  \lambda^{5/2} \sqrt{\int_{W^5\cap H^{(4)}_\lambda }g(x_1,x_2)g(x_1,x_3)g(x_2,x_4)g(x_3,x_5) \theta^5(dx_1,...,dx_5)}$.

\end{enumerate}

The following expressions are obtained by setting $g=1$:

\begin{enumerate}

\item[\rm (i')] $V^2_\lambda(1) =V^2_{1,\lambda}(1) +V^2_{2,\lambda}(1) =4 \lambda^3 \theta(W^3\cap H^{(1)}_\lambda) + 2 \lambda^2\theta^2(W^2\cap H_\lambda)$;

\item[\rm (ii')] $A_{\lambda}(1) =  \lambda^{5/2} \sqrt{\theta^5(W^5\cap H_\lambda^{(2)})}$;

\item[\rm (iii')] $B_{\lambda}(1) =  \lambda^{2} \sqrt{\theta^4 (W^4\cap H_\lambda^{(3)})}$;

\item[\rm (iv')]  $C_{\lambda}(1) = \lambda^{3/2} \sqrt{\theta^3(W^3\cap H_\lambda^{(1)})}$;

\item[\rm (v')] $D_{\lambda}(1) =  \lambda\sqrt{\theta^2(W^2\cap H_\lambda})$;

\item[\rm (vi')] $E_{\lambda}(1) =  \lambda^{5/2} \sqrt{\theta ^5(W^5\cap H^{(4)}_\lambda )}$.

\end{enumerate}

\subsection{General conditions and bounds}\label{ss:genbounds}

We start with a general estimate.

\begin{thm}[General bound for geometric graphs] 
\label{thm:bound-ABCDE}
Let the previous assumptions and notation prevail, and let $N\sim \mathscr{N}(0,1)$. Then, there exists a universal constant $C$, not depending on $\lambda$, such that, for every $\lambda >0$,
\begin{equation}\label{e:b1}
d_W(\widetilde{F}(g,W; \eta_\lambda,H_\lambda),N) \leq \frac{C_0}{V_\lambda^2(g)}\times \max  \{ A_\lambda(g), B_\lambda(g), C_\lambda(g), D_\lambda(g), E_\lambda(g)\}.
\end{equation}
If the class $\{\widetilde{F}(g,W; \eta_\lambda,H_\lambda)^4 : \lambda>0\}$ is uniformly integrable, then the RHS of {\rm (\ref{e:b1})} converges to zero, as $\lambda\to \infty$, if and only if the CLT {\rm (\ref{e:graphclt}) }takes place.
\end{thm}
\noindent{\it Proof.} 
In what follows, we write $F_\lambda = F(g,W; \eta_\lambda,H_\lambda)$ and $\widetilde{F}_\lambda = \widetilde{F}(g,W; \eta_\lambda,H_\lambda)$ to simplify the notation. Last and Penrose's formula (\ref{e:lastpenrose}) implies that the random variable $F_\lambda$ admits the following chaotic decomposition
\begin{equation}\label{e:2C}
F_\lambda = E[F_\lambda] +I_1(f_{1,\lambda}) +I_2(f_{2,\lambda}) := E[F_\lambda] +F_{1,\lambda} +F_{2,\lambda} ,
\end{equation}
where $f_{1,\lambda}(x) =2\lambda \int_Z {\bf 1}_{\{H_\lambda\cap W\times W\}} (x,y)g(x,y)\theta(dy)$ and $f_{2,\lambda}(x_1,x_2) = {\bf 1}_{H_\lambda\cap (W\times W)} (x_1,x_2)g(x_1,x_2)$. Routine computations imply then that
\begin{equation}\label{e:variances}
V^2_{i,\lambda}(g) = \var(F_{i,\lambda}), \quad i=1,2. 
\end{equation}
It follows that 
\[
\widetilde{F}_\lambda = I_1(f_1/V_{\lambda}(g)) +I_2(f_2/V_{\lambda}(g)).
\]
The upper bound (\ref{e:b1}) is now obtained by using (\ref{e:B3}), as well as the following relations, that can be proved by a standard use of the Fubini Theorem:
\begin{eqnarray*}
 4A_\lambda(g) &=&  \| f_{1,\lambda}\|^2_{L^4(\mu_\lambda)},\\
 B_\lambda(g) &=& \| f_{2,\lambda}\star_1^1 f_{2,\lambda} \|_{L^2(\mu^2_\lambda)},\\
 C_\lambda(g) &=& \| f_{2,\lambda}\star_2^1 f_{2,\lambda} \|_{L^2(\mu_\lambda)},\\
 D_\lambda(g) &=& \| f_{2,\lambda} \|^2_{L^4(\mu^2_\lambda)},\\
 2E_\lambda(g) &=& \| f_{1,\lambda}\star_1^1 f_{2,\lambda} \|_{L^2(\mu_\lambda)}.
\end{eqnarray*} 
The last assertion in the statement follows from a direct application of Theorem \ref{t:main2}.
\fin

The next two statements provide simplified bounds in case one of the two elements of the chaotic decomposition of $\widetilde{F}(g,W; \eta_\lambda,H_\lambda)$ converges to zero, as $\lambda \to \infty$. The proof (which is standard and left to the reader) uses (\ref{e:variances}) as well as the following basic estimate: if $Q,R,S$ are three random variables in $L^1(P)$, then
\[
d_W(Q+R,S) \leq d_W(R,S) +E[ |Q |].
\]

\begin{prop}[Dominating first chaos] If
\begin{equation}\label{e:firstdominates}
\frac{V_{2,\lambda}(g)}{V_{1,\lambda}(g)} \rightarrow 0,\quad \mbox{as}\quad \lambda\to\infty,
\end{equation}
then there exists a constant $C_1$, independent of $\lambda$, such that, for $\lambda$ large enough,
\begin{equation}\label{e:b2}
d_W(\widetilde{F}(g,W; \eta_\lambda,H_\lambda),N) \leq C_1\left\{\frac{V_{2,\lambda}(g)}{V_{1,\lambda}(g)} +  \frac{1}{V_{1,\lambda}^2(g)}\times \max  \{ A_\lambda(g)\}\right\}.
\end{equation}
If the class $\{\widetilde{F}(g,W; \eta_\lambda,H_\lambda)^4 : \lambda>0\}$ is uniformly integrable, then the RHS of {\rm (\ref{e:b2})} converges to zero, as $\lambda\to \infty$, if and only if the CLT {\rm (\ref{e:graphclt})} takes place.
\end{prop}

\begin{prop}[Dominating second chaos]
If
\begin{equation}\label{e:seconddominates}
\frac{V_{1,\lambda}(g)}{V_{2,\lambda}(g)} \rightarrow 0,\quad \mbox{as}\quad \lambda\to\infty,
\end{equation}
then there exists a constant $C_2$, independent of $\lambda$, such that, for $\lambda$ large enough,
\begin{equation}\label{e:b3}
d_W(\widetilde{F}(g,W; \eta_\lambda,H_\lambda),N) \leq C_2\left\{\frac{V_{1,\lambda}(g)}{V_{2,\lambda}(g)} +  \frac{1}{V_{2,\lambda}^2(g)}\times \max  \{ B_\lambda(g),C_\lambda(g),D_\lambda(g))\}\right\}.
\end{equation}
If the class $\{\widetilde{F}(g,W; \eta_\lambda,H_\lambda)^4 : \lambda>0\}$ is uniformly integrable, then the RHS of {\rm (\ref{e:b3})} converges to zero, as $\lambda\to \infty$, if and only if the CLT {\rm (\ref{e:graphclt})} takes place.
\end{prop}

\subsection{Edge counting in stationary graphs}\label{sec:stat-graphs}

For the rest of the section, we fix an integer $d\geq 1$. For every $\lambda >0$, we define the set $Q_{\lambda}=[-\frac{1}{2}\lambda^{1/d},\frac{1}{2}\lambda^{1/d}]^d$. We now specialize the framework of the previous two sections to the following setting
\begin{equation}
Z=\mathbb{R}^d,\quad  W=Q_{1}, \quad g=1,\quad  \mu_{\lambda}=\lambda \leb, 
\end{equation}
where $\leb$ is the Lebesgue measure on $\mathbb{R}^d$. We shall assume that, for every $\lambda >0$, the symmetric non-diagonal set $H_\lambda$ has the form
\[
H_\lambda = \{ (x,y) : x-y \in \fact{H}_\lambda\},
\]
for some set $\fact{H}_\lambda$ verifying 
\begin{equation}\label{e:posleb}
\leb(\fact H_{\lambda}\cap Q_{1})>0.
\end{equation}

\begin{rem}{\rm 

We insist that the novelty here (with respect to the usual setting of disk graphs -- see e.g. \cite[Chapter 3]{penrosebook} and the references therein) is that $\overline{H}_{\lambda }$ need not be bounded, allowing for arbitrarily distant points to be connected. This is especially relevant whenever $\overline{H}_{\lambda }=\alpha _{\lambda } \overline{H}_{1}$, where $\alpha _{\lambda }$ is a scaling factor and $\overline{H}_{1}$ is a fixed unbounded geometric connection rule.  Unlike in the classical literature of stochastic geometry, e.g. in stabilization theory, this allows for models with unbounded interactions, such as  between distant particles. As already recalled, our approach is further applied in \cite{LRP2}, where $U$-statistics with general stationary kernels (not only taking values $0$ or $1$), and general order $k\geq 2$, are considered.
}
\end{rem}

For every $\lambda >0$, we shall write 
\begin{eqnarray}
&& F_\lambda = F(1,Q_1;\eta_\lambda, H_\lambda)\quad \text{and} \quad \widetilde{F}_\lambda = \widetilde{F}(1,Q_1;\eta_\lambda, H_\lambda), \label{e:rename}
\end{eqnarray}
where we used the notation introduced in (\ref{e:graphstat})--(\ref{e:graphclt}). With this notation, each $\frac12F_\lambda$ is a stationary $U$-statistic (see Definition \ref{d:stationary}), counting the number of edges in the stationary random graph based on $H_\lambda$ (see Definition \ref{d:randomg}). The chaotic decomposition of $F_\lambda$ is written
\[
F_\lambda = E[F_\lambda] +F_{1,\lambda}+ F_{2,\lambda},
\]
where we have adopted the same notation as in (\ref{e:2C}). Since $g=1$, Problem \ref{p:graphs} becomes the following: {\sl characterize all collections of sets $\{\fact H_\lambda\}$ such that the CLT {\rm (\ref{e:graphclt})} takes place, and assess the rate of convergence in the Wasserstein distance.}

\begin{rem}
\label{rem:rescaling}
{\rm
For every $\lambda >0$, one has the equality in law

\begin{equation}
\label{eq:def-F-stat-graphs}
F_{\lambda}\stackrel{\rm Law}{=}\sum_{x,y\in \eta\cap Q_{\lambda}, x\neq y}\1_{x-y\in \fact G_{\lambda}}, \lambda>0,
\end{equation}
where $\eta$ is a random Poisson measure with Lebesgue intensity, and $\fact G_{\lambda}$ is a measurable subset of $\mathbb{R}^d$ defined by the relation 
\begin{equation}\label{e:glambda}
\fact H_{\lambda}=\lambda^{-1/d}\fact G_{\lambda},
\end{equation}
so that
\begin{equation}
\label{eq:H-nonnegligible}
\leb(\fact G_{\lambda}\cap Q_{\lambda})>0.
\end{equation} }
\end{rem}

\begin{rem}[Asymptotic equivalence notation]{\rm
Given two mappings $\lambda \mapsto \gamma_{\lambda},\,  \lambda \mapsto \delta_{\lambda}$, we write $\gamma_{\lambda}\asymp \delta_{\lambda}$ if there are two positive constants $C,C'>0$ such that $C \gamma_{\lambda} \leq \delta_{\lambda} \leq C' \gamma_{\lambda}$ for $\lambda$ sufficiently large. We write $\gamma_{\lambda}\sim\delta_{\lambda}$ if $\delta_{\lambda}>0$ for $\lambda$ sufficiently large and $\alpha_{\lambda}/\delta_{\lambda}\to 1$.}
\end{rem}

One of the main points developed in the present section is that the asymptotic Gaussianity of the class $\{\widetilde F_\lambda\}$ results can be effectively studied by using the \emph{occupation coefficient} of $\fact H_{\lambda}$, defined as
\begin{equation}
\label{eq:def-psi}
\psi(\lambda):=\leb(\fact H_{\lambda}\cap Q_{1})=\frac{\leb(\fact G_{\lambda}\cap Q_{\lambda})}{\leb(Q_{\lambda})},\quad \lambda>0.
\end{equation}
We also write $\hat W=W-W=Q_{2^d}$ and $\check W=Q_{1/2^d}$ (this is the largest set such that $\check W-\check W\subseteq W$) and define, for $\lambda>0$, 
\begin{align}
\label{eq:def-check-psi}
\check \psi(\lambda):=\leb(\fact H_{\lambda}\cap \check W)=\frac{\leb(\fact G_{\lambda}\cap Q_{\lambda/2^d})}{\leb(Q_{\lambda})},\quad \lambda>0,
\\
\label{eq:def-hat-psi}
\hat \psi(\lambda):=\leb(\fact H_{\lambda}\cap \hat W)=\frac{\leb(\fact G_{\lambda}\cap Q_{2^d \lambda})}{\leb(Q_{\lambda})}, \quad \lambda>0.
\end{align}

In order to obtain necessary and sufficient conditions for asymptotic normality, we will often work under the additional assumption that 
\begin{equation}
\label{eq:psi-regular}
\check \psi(\lambda)\asymp\hat  \psi(\lambda).
\end{equation}   
In this case, one has trivially that $\check \psi(\lambda)\asymp\hat  \psi(\lambda)\asymp \psi(\lambda)$,
and the value of $\psi$ is only relevant up to a fixed multiplicative constant.

\begin{rem}[O-regularity]\label{r:oreg}{\rm
Assume that the geometric rule defined by $\fact G_{\lambda}$ does not depend on $\lambda$, i.e.: $\fact G_{\lambda}=\fact G$ for some fixed measurable set $\fact G$, in such a way that each set $\fact H_\lambda$ is obtained by rescaling $\fact G$ by a factor $\lambda^{-1/d}$. Then, condition (\ref{eq:psi-regular}) is implied by the following stronger assumption: $\psi(a \lambda)\asymp \psi(\lambda)$ for every $a>0$. In the terminology of \cite[Section 2.2]{Ency-27}, this is equivalent to saying that $\psi$ is \emph{O-regular}. }
\end{rem}

In view of using the bounds appearing in Theorem \ref{thm:bound-ABCDE}, we have the following crucial estimates:

\begin{thm}\label{t:esti}Let the previous notation and assumption prevail,
set $V_{i,\lambda}(1) =V_{i,\lambda}$, $i=1,2$, as well as  $A_{\lambda}=A_{\lambda}(1), B_{\lambda}=B_{\lambda}(1),\dots, E_{\lambda}=E_{\lambda}(1)$  -- see Section {\rm \ref{ss:frame}}.
The following estimates are in order for every fixed $\lambda >0$:
\begin{eqnarray*}
&& 2^{-d}\check\psi({ \lambda}) \leq \lambda^{-2}V_{2,\lambda}^2=2\lambda^{-2}D_{\lambda}^2 \leq  2^d \hat\psi( \lambda),\\
&&2^{-d}\check\psi({ \lambda})^2 \leq \lambda^{-3}V_{1,\lambda}^2=\lambda^{-3}C_{\lambda}^2 \leq  2^d \hat\psi( \lambda)^2,\\
&&2^{-d}\check\psi({ \lambda})^4 \leq \lambda^{-5}A_{\lambda}^2 \leq  2^d \hat\psi( \lambda)^4,\\
&&  \lambda^{-4}B_{\lambda}^2 \leq  2^d\hat\psi( \lambda)^3,\\
&& 2^{-d}\check\psi( \lambda)^4  \leq \lambda^{-5}E_{\lambda}^2 \leq  2^d\hat\psi( \lambda)^4.
\end{eqnarray*}

 \end{thm}

\noindent {\it Proof.}
We introduce the   changes of variables denoted by $\varphi^{(i)}, i=0,1,2,3,4$, where

\begin{align*}
\varphi^{(0)}(x_{1},x_{2})&=(x_{1},u=x_{1}-x_{2}),\\
\varphi^{(1)}(x_{1},x_{2},x_{3})&=(x_{1},u=x_{1}-x_{2},v=x_{1}-x_{3}),\\
\varphi^{(2)}(x_{1},x_{2},x_{3},x_{4},x_{5})&=(x_{1},u=x_{1}-x_{2},v=x_{1}-x_{3},w=x_{1}-x_{4},z=x_{1}-x_{5}),\\
\varphi^{(3)}(x_{1},x_{2},x_{3},x_{4})&=(x_{1},u=x_{1}-x_{2},v=x_{2}-x_{3},w=x_{3}-x_{4}),\\
\varphi^{(4)}(x_{1},x_{2},x_{3},x_{4},x_{5})&=(x_{1},u=x_{1}-x_{2},v=x_{2}-x_{3},w=x_{3}-x_{4},z=x_{4}-x_{5}).
\end{align*}
Using the notation introduced in Definition \ref{def:ABCDE} we have
\begin{align*}
\leb(W^2\cap H_{\lambda})&=\int_{W^2}\1_{x_{1}-x_{2}\in \fact H_\lambda}dx_{1}dx_{2}=\int_{\varphi^{(0)}(W^2)}\1_{u\in \fact H_\lambda}dx_{1}du,\\
\leb(W^3\cap H_{\lambda}^{(1)})&=\int_{W^3}\1_{x_{1}-x_{2}\in \fact H_\lambda}\1_{x_{1}-x_{3}\in \fact H_{\lambda}}dx_{1}dx_{2}dx_{3}=\int_{\varphi^{(1)}(W^3)}\1_{u\in \fact H_\lambda}\1_{v \in \fact H_\lambda}dx_{1}dudv,\\
\leb(W^5\cap H_{\lambda}^{(2)})&=\int_{W^5}\1_{x_{1}-x_{2}\in \fact H_\lambda\1}\1_{x_{1}-x_{3}\in \fact H_\lambda}\1_{x_{1}-x_{4}\in \fact H_\lambda}\1_{x_{1}-x_{5}\in \fact H_\lambda}dx_{1}dx_{2}dx_{3}dx_{4}dx_{5}\\
&=\int_{\varphi^{(2)}(W^5)}\1_{u\in \fact H_\lambda}\1_{v \in \fact H_\lambda}\1_{w\in \fact H_\lambda}\1_{z \in \fact H_\lambda}dx_{1}dudvdwdz,\\
\leb(W^4\cap H_{\lambda}^{(3)})&=\int_{W^4}\1_{x_{1}-x_{2}\in \fact H_\lambda}\1_{x_{2}-x_{3}\in \fact H_\lambda}\1_{x_{3}-x_{4}\in \fact H_\lambda}\1_{x_{4}-x_{1}\in \fact H_\lambda}dx_{1}dx_{2}dx_{3}dx_{4}\\
&=\int_{\varphi^{(3)}(W^4)}\1_{u\in \fact H_\lambda}\1_{v \in \fact H_\lambda}\1_{w\in \fact H_\lambda}\1_{u+v+w \in \fact H_\lambda}dx_{1}dudvdw,\\
\leb(W^5\cap H_{\lambda}^{(4)})&=\int_{W^5}\1_{x_{1}-x_{2}\in \fact H_\lambda}\1_{x_{2}-x_{3}\in \fact H_\lambda}\1_{x_{3}-x_{4}\in \fact H_\lambda}\1_{x_{4}-x_{5}\in \fact H_\lambda}dx_{1}dx_{2}dx_{3}dx_{4}dx_{5}\\
&=\int_{\varphi^{(4)}(W^5)}\1_{u\in \fact H_\lambda}\1_{v \in \fact H_\lambda}\1_{w\in \fact H_\lambda}\1_{z \in \fact H_\lambda}dx_{1}dudv.\\
\end{align*}

Using the inclusions 
\begin{align*}
\check W^2& \subseteq \varphi^{(0)}(W^2) \subseteq \hat W^2,\\
\check W^3& \subseteq \varphi^{(1)}(W^3) \subseteq \hat W^3,\\
\check W^5& \subseteq\varphi^{(2)}(W^5) \subseteq \hat W^5,\\
& \varphi^{(3)}(W^4) \subseteq \hat W^4,\\
\check W^5 & \subseteq \varphi^{(4)}(W^5) \subseteq \hat W^5,\\
\end{align*}
and

\begin{equation*}
 \1_{u\in  \fact H_\lambda}\1_{v \in  \fact H_\lambda}\1_{w\in  \fact H_\lambda}\1_{u+v+w \in  \fact H_\lambda} \leq \1_{u\in  \fact H_\lambda}\1_{v \in  \fact H_\lambda}\1_{w\in  \fact H_\lambda},\, u,v,w \in \mathbb{R}^d,
\end{equation*}

we have 

\begin{align*}
\leb^2(\check W\times ( \fact H_\lambda\cap \check W))=2^{-d}\check\psi({ \lambda})& \leq \lambda^{-2}V_{2,\lambda}^2 \leq \leb^2(\hat W\times ( \fact H_\lambda\cap \hat W))= 2^d \hat\psi( \lambda)\\
\leb^3(\check W\times ( \fact H_\lambda\cap \check W)^2)=2^{-d}\check\psi({ \lambda})^2& \leq  \lambda^{-3}C_{\lambda}^2 \leq \leb^3(\hat W\times ( \fact H_\lambda\cap \hat W)^2)= 2^d \hat\psi( \lambda)^2\\
\leb^5(\check W\times ( \fact H_\lambda\cap \check W)^4)=2^{-d}\check\psi({ \lambda})^4& \leq  \lambda^{-5}A_{\lambda}^2 \leq \leb^5(\hat W\times ( \fact H_\lambda\cap \hat W)^4)= 2^d \hat\psi( \lambda)^4\\
 &  \lambda^{-4} B_{\lambda}^2 \leq \leb^4(\hat W\times ( \fact H_\lambda\cap \hat W)^3)= 2^d \hat\psi( \lambda)^3\\
\leb^5(\check W\times ( \fact H_\lambda\cap \check W)^4)=2^{-d}\check\psi({ \lambda})^4 &\leq \lambda^{-5} E_{\lambda}^2 \leq \leb^5(\hat W\times ( \fact H_\lambda\cap \hat W)^4)= 2^d \hat\psi( \lambda)^4\\
\end{align*}
and the result follows.

\qed

The next statement provides one of the main results of this section: it gives an exhaustive characterization of the asymptotic behavior of $F_\lambda$, whenever (\ref{eq:psi-regular}) is in order. In order to allow for a comparison with the existing literature, we classify the asymptotic behavior of $F_\lambda$ according to four regimes, denoted by {\bf (R1)}}--{\bf (R4)}. Such a classification is based on the proportion  $\psi (\lambda )$ of space occupied by $\overline{H}_{\lambda }$ in the observation window, determining the  influence area of a given point of the Poisson measure. This coefficient has to be compared with  $\lambda ^{{-1}}$, which corresponds to the total window measure divided by the mean number of points. The four regimes are the following:
\begin{description}
\item{{\bf(R1)}} $\lambda  \psi (\lambda )\to 0$  and $\lambda \sqrt{\psi (\lambda )}\to \infty $;

\item {{\bf (R2)}}  $\lambda  \psi (\lambda )\to \infty $;

\item{{\bf (R3)}} \emph{(Thermodynamic regime)}  $\lambda \psi (\lambda )\asymp 1$;
 
\item{{\bf(R4)}}  The mapping $ \lambda\mapsto \lambda \sqrt{\psi (\lambda )}$ is bounded.

\end{description}

The thermodynamic regime corresponds (after rescaling) to the usual models where the geometry of the interactions does not change as the window of observation grows to the whole space (see Remark \ref{rem:rescaling}). We will see in Section \ref{ss:exex} that, when specialized to Poissonized disk graphs, our asymptotic approximations and variance estimates  concur with those obtained in \cite[Chapter 3]{penrosebook}. Under regimes {\bf (R2)} and {\bf (R3)}, there is asymptotic normality with convergence at speed $\lambda ^{-1/2}$ in the Wasserstein distance. Under {\bf (R1)} the  convergence to the normal law is slower, and under {\bf (R4)} the asymptotic normality is lost:  for any converging subsequence the limit is either Poisson or zero.

\begin{rem}{\rm One interesting contribution of Theorem \ref{thm:CLT-stat-graphs} appears at the end of Point (ii), where it is stated that, under the thermodynamic regime, both chaotic projections of the random variable $\tilde{F}_\lambda$ contributee to the limit and satisfy a joint CLT. This kind of phenomenon is an example of the ``fine Gaussian fluctuations'' appearing in the title of the paper.

}
\end{rem}

\begin{thm}
\label{thm:CLT-stat-graphs}
Let $\{\fact H_{\lambda} : \lambda>0\}$ be a family of subsets of $\mathbb{R}^d$ satisfying {\rm (\ref{e:posleb})} and let $\check\psi,\, \hat \psi$ be defined according to {\rm  (\ref{eq:def-check-psi})--(\ref{eq:def-hat-psi})}.
Assume in addition that {\rm (\ref{eq:psi-regular})} is satisfied, and consider a random variable $N\sim \mathscr{N}(0,1)$.
The quantities introduced in Section {\rm \ref{ss:frame}} satisfy the following relations: there exist constants $0<k<K<\infty$, independent of $\lambda$, such that
\begin{eqnarray*}
 \frac{V_{1,\lambda}^2}{V_{2,\lambda}^2}\asymp\lambda \psi(\lambda),\,\, && \frac{1}{V_{1,\lambda}^2}A_{\lambda} \asymp  \lambda^{-1/2}, \quad \text{and}\\
k\max(\lambda^{-1/2},
%\textcolor{red}{  \lambda^{-1/2},},
(\lambda^2\psi(\lambda))^{-1/2},\lambda^{1/2}\psi(\lambda))&\leq&\frac{1}{V_{2,\lambda}^2}\max(B_{\lambda},C_{\lambda},D_{\lambda},E_{\lambda}) \\
&\leq& K \max(\sqrt{ \psi(\lambda)}, \lambda^{-1/2}, (\lambda^2\psi(\lambda))^{-1/2},\lambda^{1/2}\psi(\lambda))).
\end{eqnarray*}
Furthermore, one can  choose $K$ in such a way that the following properties {\rm (i)--(iii)} are verified.
\begin{description}
\item[\rm (i)] (Regime {\bf (R2)}) If $\lambda \psi(\lambda) \to \infty$, the first chaos projection $ F_{1,\lambda}$ dominates and 
\begin{equation*}
\var(F_{1,\lambda})\sim \var(F_{\lambda})\asymp \lambda^3 \psi(\lambda)^2 \to \infty,
\end{equation*} and 
\begin{align*}
%d_{W}(\tilde F_{1,\lambda},N) & \leq C\lambda^{-1/2},\\
d_{W}(\tilde F_{\lambda},N)  &  \leq K\, \lambda^{-1/2}.
\end{align*}

\item[\rm (ii)] (Regime {\bf (R3)}) If $\lambda \psi(\lambda) \asymp 1$,  
\begin{equation*}
\var(F_{\lambda}) \asymp \var(F_{1,\lambda})\asymp \var(F_{2,\lambda})\asymp \lambda,
\end{equation*} and
\begin{equation*}
d_{W}(\tilde F_{\lambda},N) \leq K\, \lambda^{-1/2}.
%,\quad d_{W}(\tilde F_{1,\lambda},N) \leq C \lambda^{-1/2},\quad d_{W}(\tilde F_{2,\lambda},N) \leq C \lambda^{-1/2},
\end{equation*}
In this case one has also that, as $\lambda\to\infty$, the pair 
\[
\left(\frac{F_{1,\lambda}}{\var(F_{1,\lambda})^{1/2}},\,\, \frac {F_{2,\lambda}}{\var(F_{2,\lambda})^{1/2}} \right),
\]
converges in distribution to a two-dimensional Gaussian vector $(N_1,N_2)$, such that $N_i \sim \mathscr{N}(0,1)$ and $N_1,N_2$ are independent. 
\item[\rm (iii)] (Regimes {\bf (R1)} and {\bf (R4)}) If $\lambda \psi(\lambda)\to 0$, then the second chaos projection $F_{2,\lambda}$ dominates,
\begin{equation*}
\var(F_{\lambda})\sim\var( F_{2,\lambda})\asymp \lambda^2 \psi(\lambda) ,
\end{equation*}
and
\begin{equation*}
d_{W}(\tilde F_{\lambda}, N) \leq \frac{ K}{\lambda\sqrt{\psi(\lambda)}}.
\end{equation*}
Moreover, one has that $\tilde F_{\lambda}$ converges in distribution to $N$ if and only if $\lambda^2 \psi(\lambda)\to \infty$ (or, equivalently, $\var(F_{\lambda})\to \infty$).
\end{description}

\end{thm}

%\begin{rem}{\rm It is plausible that the case where the mapping $\lambda \mapsto \lambda^2\psi(\lambda)$ is {\it not} bounded away from zero cannot be studied by means of the techniques developed in this paper. Indeed, the estimate (\ref{e:4mb}) implies that, for large $\lambda$, the quantity $(\lambda^2\psi(\lambda))^{-1}$ is a rough lower bound for the difference $E[\tilde F_\lambda^4] -3$. We therefore conjecture that, in this case, the sequence $\tilde F_\lambda$ converges in distribution to some random variable $X$ such that $E[X^4] = \infty$. We leave this point open for future research.
%}
%\end{rem}

\noindent{\it Proof of Theorem {\rm \ref{thm:CLT-stat-graphs}}.} Theorem \ref{thm:bound-ABCDE} and assumption (\ref{eq:psi-regular}) yield
\begin{align*}
V_{2,\lambda}^2=2D_{\lambda}^2& \asymp \lambda^{2}\psi( \lambda)\\
V_{1,\lambda}^2=C_{\lambda}^2 &\asymp \lambda^3\psi(\lambda)^2\\
 A_{\lambda}^2 &\asymp \lambda^5 \psi( \lambda)^4\\
 B_{\lambda}^2& \leq 2^d \lambda^4 \psi( \lambda)^3\\
E_{\lambda}^2 &\asymp \lambda^5\psi( \lambda)^4.
\end{align*}
(i) We assume $\lambda \psi(\lambda)\to \infty$. Applying  (\ref{e:b1}), one deduces that 
\begin{align*}
d_{W}(\tilde F,N) & \leq \frac{1}{V_{1,\lambda}^2}\max(A_{\lambda},B_{\lambda},C_{\lambda},D_{\lambda},E_{\lambda})\\
& \leq \frac{K\sqrt{\lambda}}{\lambda^3 \psi(\lambda)^2}\max(\lambda^{2}\psi(\lambda)^2, \lambda^{3/2} \psi(\lambda)^{3/2}, \lambda\psi(\lambda),\lambda^{1/2} \psi(\lambda)^{1/2})\\
& \leq \frac{K}{\sqrt{\lambda}}
\end{align*}
for $\lambda$ large enough.  Notice that (\ref{e:b2}) does not yield a better bound because $\psi(\lambda) \leq 1$.

 \noindent (ii) Applying again (\ref{e:b1}), the conclusion is deduced from Point (i), because $\lambda\psi(\lambda)\geq c$ for some constant $c>0$ and for $\lambda$ large enough. The last statement at Point (ii) follows from an application of Proposition \ref{p:sumtosingle}.
 
\noindent (iii) Using (\ref{e:b1}) again yields 
\begin{align*}
d_{W}(\tilde F, N) & \leq \frac{1}{\lambda^2 \psi(\lambda)} \max(\lambda^{5/2}\psi(\lambda)^2, \lambda^2 \psi(\lambda)^{3/2}, \lambda^{3/2}\psi,\lambda \psi(\lambda)^{1/2})\\
& \leq \frac{K}{\lambda\sqrt{\psi(\lambda)}}
\end{align*}
 for $\lambda$ large enough, because $\lambda \psi(\lambda)\to 0$. To conclude the proof, we have to show that, if $\lambda^2 \psi(\lambda)$ does not diverge to infinity, then $\tilde F_{\lambda}$ does not converge in distribution to $N$. To prove this negative result, one could apply the product formula (\ref{e:product}) to prove that, whenever $\lambda^2 \psi(\lambda)$ is not diverging to infinity and is bounded away from zero, there exists a sequence $\lambda_n$, $n\geq 1$, such that $\lambda_n \to \infty$ and $\sup_n E[\tilde F_{\lambda_n}^6]<\infty$, so that the desired conclusion is deduced from the last part of Theorem \ref{thm:bound-ABCDE} (the case when $\lambda^2 \psi(\lambda)$ is not bounded away from zero can be dealt with by a direct argument). However, the statement of the forthcoming Theorem \ref{t:poisson} is much stronger, and it is therefore not necessary to spell out the details of these computations.
\qed

\begin{cor}
Assume that the geometric rule defined by $\fact G_{\lambda}$ (see {\rm (\ref{e:glambda})}) does not depend on $\lambda$, in such a way that $\fact G_{\lambda}=\fact G$ for some fixed measurable set $\fact G$. Assuming {\rm (\ref{eq:psi-regular})} (see Remark {\rm \ref{r:oreg}}), one has that $\tilde F_{\lambda}$ converges in distribution to $N\sim\mathscr{N}(0,1)$, with a rate at most of the order $\lambda^{-1/2}$ with respect to $d_W$.  
\end{cor}

\noindent{\it Proof}. We are in one of the following situations:
\begin{enumerate}
\item If $\fact G$ has finite positive Lebesgue measure, $\check \psi(\lambda)\asymp \hat \psi(\lambda)\asymp \leb(\fact G)/\leb(Q_{\lambda})\asymp\lambda^{-1}$. It corresponds to the case (ii) in Theorem \ref{thm:CLT-stat-graphs}, meaning the two chaoses codominate. It follows that $\tilde F_{\lambda}$ converges to the normal law with a rate at most of the order of $\lambda^{-1/2}$ with respect to $d_{W}$.
\item If $\fact G$ does not have finite measure, $\lambda \psi(\lambda)\asymp \leb(\fact G\cap Q_{\lambda})\to \infty$ and we are in the situation of Point (i) of Theorem \ref{thm:CLT-stat-graphs}, that is: the first chaos dominates. We therefore deduce that 
\begin{equation*}
d_{W}(\tilde F_{\lambda},N) \leq K \lambda^{-1/2},
\end{equation*}
for some $K>0$, and the conclusion follows.
\end{enumerate}

\qed

%Besides $\lambda\psi(\lambda)\asymp 1$, the second cutting point is  $\psi(\lambda) \asymp  \lambda^{-2}$. It corresponds for instance to $\fact G_{\lambda}=\lambda^{-1/d}\fact G$ for some measurable $\fact G$ with positive finite Lebesgue measure. Within a standard Poisson measure on $\mathbb{R}^d$, the random integer number  $F_{\lambda}$ of pair of points $(x,y)$ in $Q_{\lambda}$ such that $y\in x+\lambda^{-1/d}\fact H$ has a bounded variance (or at least ${\tilde F_{\lambda_{n}};\,n \geq1}$ has a bounded variance for some sequence $\lambda_{n}\uparrow \infty$) and does not follow a CLT.
%%variable with bounded variance, whence $\tilde F_{\lambda}$ cannot be asymptotically normal. Every situation where there are even less points cannot be asymptotically normal either.
%
%Many things can be said about $\tilde F_{\lambda},\tilde F_{1,\lambda},\tilde F_{2,\lambda}$ if (\ref{eq:psi-regular}) is not satisfied by mimicking the arguments of Th. \ref{thm:CLT-stat-graphs}. We simply give the following result that seems to be the most important fact.

As announced, we shall now deal more thoroughly with the case where $\lambda^2 \psi(\lambda)$ does not diverge to infinity. In the proof of the next statement we shall use the following notation: if $X$ is a random variable with finite moments of every order, then we write $\{\chi_m(X) : m\geq 1\}$ to indicate the sequence of its cumulants (see \cite[Chapter 3]{PeTa} for an introduction to this concept). For instance, $\chi_1(X) = E[X]$, $\chi_2(X) = \var(X)$, and so on.

\begin{rem}{\rm The proof of Theorem \ref{t:poisson} provided below is based on diagram formulae and the method of moments and cumulants. An alternate proof could be deduced from the classic results by Silverman and Brown \cite{SiBr}, combined with a Poissonization argument. Another proof of this result, complete with explicit bounds in the total variation distance, appears in \cite{Pec2011}. The proof provided below has the merit of illustrating an application of diagram formulae (that are typically used to deduce CLTs) to a non-central result.
}
\end{rem}

\begin{thm}[Poisson approximations]\label{t:poisson} Let the assumptions and notation of Theorem {\rm \ref{thm:CLT-stat-graphs}} prevail, and assume in addition that $\lambda\psi(\lambda)\to 0$ and the mapping $\lambda \mapsto \lambda^2\psi(\lambda)$ does not diverge to infinity, as $\lambda \to \infty$. Then, there exists a sequence $\{\lambda_n\}$ such that $\lim_{n\to\infty}\lambda_n = \infty$ and the following properties {\rm (i)-(iii)} hold.
\begin{itemize}

\item[\rm (i)] $\lim_{n\to\infty}\var(F_{\lambda_{n}}) =\lim_{n\to\infty}\var(F_{2,\lambda_{n}})  = 2c$, where $c$ is a nonnegative constant possibly depending on the sequence $\{\lambda_n\}$, and we have used the notation {\rm (\ref{e:2C})}.

\item[\rm (ii)] If $c=0$, then $E[ |\tilde F_{\lambda_{n}}|] \to 0$, as $n\to\infty$.

\item[\rm (iii)] If $c>0$, then $ F_{2,\lambda_{n}}$, and therefore $F_{\lambda_{n}} - E[F_{\lambda_{n}}]$, converges in distribution to $X =2 P(c/2)$, as $n\to\infty$, where $P(c/2)$ indicates a centered Poisson random variable with parameter $c/2$. 
\end{itemize}

In particular, the family $\{\widetilde F_\lambda\}$ does not verify a CLT as $\lambda \to \infty$.

\end{thm}
{\it Proof.} Since $\var(F_\lambda)\sim \var(F_{2,\lambda}) \asymp \lambda^2\psi(\lambda)$, Point (i) is immediately deduced from the Bolzano-Weierstrass theorem. Point (ii) follows from a direct application of Campbell's Theorem (see \cite[Theorem 3.1.3]{SchnWeil}), yielding that, as $n \to \infty$,
\[
\frac{E[F_{\lambda_{n}}]}{\var(F_{\lambda_{n}})^{1/2}}\asymp \lambda_{n} \sqrt{\psi(\lambda_{n})}\to 0.
\]
We shall prove Point (iii) by using the method of cumulants. First of all, we observe that since $\lambda_{n}^2\psi(\lambda_{n})$ is bounded and bounded away from zero, one has that $V^2_{1,\lambda_{n}} \asymp \lambda_n^3\psi(\lambda_n)^2 \to 0$, that is: as $n\to\infty$, the limits of $ F_{2,\lambda_{n}}$and $F_{\lambda_{n}} - E[F_{\lambda_{n}}]$ coincide. We recall that the law of the random variable $X = 2P(c/2)$ is determined by its moments or, equivalently, by its cumulants (see e.g. \cite[pp. 42-43]{PeTa}). Standard computations imply that $\chi_1(X) = 0$ and, for every $m\geq 2$, $\chi_m(X) = 2^{m-1} c$. We are therefore left to show that, for every $m\geq 3$,
\begin{equation}\label{e:inproof}
\chi_m(I_2(f_{2,\lambda_n})) \longrightarrow 2^{m-1} c,
\end{equation}
where $f_{2,\lambda}$ is the kernel appearing in formula (\ref{e:2C}). Our proof of (\ref{e:inproof}) is based on the use of the so-called `diagram formulae' for multiple Poisson integrals as stated in \cite[Corollary 7.4.1]{PeTa} (observe that this statement only deals with simple kernels: however, one can immediately extend this result to our framework by using e.g. Surgailis \cite[Theorem 3.1]{surgailis84}). Fix $m\geq 3$, and introduce the following notation and definitions.  We shall write $[2m] = \{1,2,...,2m\}$ and indicate by $\pi_0$ the partition of $[2m]$ given by $\pi_0 = \{\{1,2\},\{3,4\},...,\{2m-3,2m-2\},\{2m-1, 2m\}\}$. Given a partition $\pi$ of $[2m]$, we write $|\pi|$ to indicate the number of blocks of $\pi$. The function $(f_{2,\lambda})_\pi$, in $|\pi|$ variables, is obtained as follows: (1) consider the function 
\[
\Phi(x_1,....,x_{2m}) = f_{2,\lambda}(x_1,x_2)\times f_{2,\lambda}(x_3,x_4)\times \cdots \times f_{2,\lambda}(x_{2m-1},x_{2m}),
\]
given by the juxtaposition of $m$ copies of $f_{2,\lambda}$, and (2) identify two variables $x_i , x_j$ in the argument of $\Phi$ if and only if $i$ and $j$ are in the same block of $\pi$. According to \cite[Corollary 7.4.1]{PeTa}, one has therefore that
\[
\chi_m(I_2(f_{2,\lambda_n})) = \sum_{\pi \in \mathcal{M}_m} \lambda_n^{|\pi|} \int_{(\R^d)^{|\pi|}} (f_{2,\lambda_n})_\pi \,d\leb^{|\pi|},
\]
where the symbol $\mathcal{M}_m$ stands for the class of those partitions $\pi$ of $[2m]$ satisfying the following properties: (a) every block of $\pi$ contains at least two elements, (b) given any two blocks $b_0\in \pi_0$ and $b_1 \in \pi$, the intersection $b_0\cap b_1$ contains at most one element,  and (c) the diagram $\Gamma(\pi_0, \pi) $, as defined in \cite[Section 4.1]{PeTa}, is connected in the sense of \cite[p. 47]{PeTa}. There are exactly $2^{m-1}$ partitions $\pi\in \mathcal{M}_m$ such that $|\pi| =2$, and for any such partition one has that
\[
\lambda_n^{|\pi|} \int_{(\R^d)^{|\pi|}} (f_{2,\lambda_n})_\pi \,d\leb^{|\pi|} = \|f_{2, \lambda_n}\|^2_{L^2((\lambda_n\leb)^2)} \longrightarrow c.
\]
On the other hand, if $\pi\in \mathcal{M}_m$ and $|\pi|\geq 3$, a change of variables analogous to the ones defined in the proof of Theorem \ref{t:esti} yields that, for some constant $C$ independent of $n$,
\[
0\leq \lambda_n^{|\pi|} \int_{(\R^d)^{|\pi|}} (f_{2,\lambda_n})_\pi \,d\leb^{|\pi|} \leq C\, \lambda_n^{|\pi|} \psi(\lambda_n)^{|\pi|-1} =C\, \lambda_n^2 \psi(\lambda_n) \times (\lambda_n \psi(\lambda_n))^{|\pi| -2} \to 0,
\]
thus concluding the proof.
\qed

\medskip

%Assume that the constant $c$ in Point (i) of Theorem \ref{t:poisson} is strictly positive. Since the variances of the random variables $\{ F_{2,\lambda_{n(k)}}\}$ converge, one natural question is the following: {\sl does there exist a random variable $X$, inside the second Wiener chaos of $\eta $ and distributed as $2P(c/2)$, such that $ F_{2,\lambda_{n(k)}}$ converges to $X$ in $L^2$?} The next general result shows that the answer is indeed negative.
%
%\begin{prop} Let Let $H\in Z\times Z$ be a symmetric set
%
%
%\end{prop}

\subsubsection{Two examples}\label{ss:exex}

We now present some explicit examples. The notation of Section \ref{sec:stat-graphs} will prevail throughout this section.

\begin{example}[Disk graph]
\label{ex:disk-graph}
{\rm
Assume 
\begin{equation*}
\fact G_{\lambda}=B(0,r_{\lambda})
\end{equation*}
for some $r_{\lambda}>0$, meaning that two points of $\eta$ in $Q_{\lambda}$ are connected whenever their distance is smaller than $r_{\lambda}$. It yields $\psi(\lambda)\asymp r_{\lambda}^d/\lambda$ (it is easy to verify that (\ref{eq:psi-regular}) is satisfied). Then $\tilde F_{\lambda}$ is asymptotically normal iff $\lambda  r_{\lambda}^d\to \infty$, and 
\begin{equation*}
d_{W}(\tilde F_{\lambda}, N) \leq C \lambda^{-1/2}\max(1,r_{\lambda}^{-d/2}).
\end{equation*}

\noindent According to the classification based on the four regimes {\bf (R1)}--{\bf (R4)}, the above result yields the following exhaustive description of the asymptotic behavior of $F_\lambda$ (note how we are able to distinguish the contribution of each chaotic projection) :
\begin{itemize}

\item[{\bf (R1)}] If $r_\lambda \to 0$ and $\lambda r_\lambda^d\to \infty$, then ${\rm Var} (F_\lambda) \asymp \lambda r_\lambda^d$, $\tilde{F}_\lambda$ satisfies a CLT with an upper bound of the order of $(\lambda r_\lambda^d)^{-1/2}$ on the Wasserstein distance, and the projection of $\tilde{F}_\lambda$ on the second Wiener chaos dominates in the limit.

\item[{\bf (R2)}] If $r_\lambda \to \infty$ , then ${\rm Var} (F_\lambda) \asymp \lambda(r_\lambda^d)^2$, $\tilde{F}_\lambda$ satisfies a CLT with an upper bound of the order of $\lambda^{-1/2}$ on the Wasserstein distance, and the projection of $\tilde{F}_\lambda$ on the first Wiener chaos dominates.

\item[{\bf (R3)}] If $r_\lambda \asymp 1$, then ${\rm Var} (F_\lambda) \asymp \lambda$, $\tilde{F}_\lambda$ satisfies a CLT with an upper bound of the order of $\lambda^{-1/2}$ on the Wasserstein distance, and the projections of $\tilde{F}_\lambda$ on the first and second Wiener chaos both contribute to the limit and satisfy a joint CLT.

\item[{\bf (R4)}] If $\lambda r^d_\lambda$ converges to a finite constant $c \geq 0$, then $F_\lambda$ converges  either to zero (if $c=0$) or to a multiple of a Poisson random variable (if $c>0$).

\end{itemize}

}
\end{example}

\begin{rem}{\rm Explicit estimates on the variances and on the rates of convergence to normal in the Wasserstein distance for edge counting statistics are also provided in \cite[Theorem 6.3]{lesmathias}.  Their variance estimates are of the order $$\max(\lambda^3 \delta_{\lambda}^{2d}, \lambda^2 \delta^d)=\max( \lambda^3 \psi(\lambda)^2, \lambda^2 \psi(\lambda) ),$$ corresponding exactly to our findings ($\delta_{\lambda} =\lambda^{-1/d}r_{\lambda}=\psi(\lambda)^{1/d}$ in their notation). They obtained convergence to the normal law if $\delta_{\lambda}\to 0$ and $\lambda^{4/3} \delta_{\lambda}^d\to \infty$, whereas Example \ref{ex:disk-graph} yields a CLT whenever $\lambda^2 \psi(\lambda)\to \infty$, dealing with the cases  $\limsup_{\lambda }\delta_{\lambda}>0$  and    $\lambda^{-2/d}=o(  \delta_{\lambda}), \delta _{\lambda } \leq  \lambda^{-4/(3d)}$, and gives a negative answer if $\delta_{\lambda} \leq C \lambda^{-2/d}$ for some $C>0$.}
\end{rem}

\begin{rem}{\rm For every fixed $\lambda$, the $U$-statistic $\frac12 F_\lambda$ has the same law as the random variable counting the number of edges in a disk graph, with radius $\delta_{\lambda} = \lambda^{-1/d}r_\lambda$, based on random points of the form $\{Y_1,...,Y_{N(\lambda)}\}$, where $\{Y_i\}$ indicates a collection of i.i.d. random variables uniformly distributed on $Q_1 = [-\frac12, \frac12]^d$, and $N(\lambda)$ is an independent Poisson random variable with parameter $\lambda$. As such, each $\frac12 F_\lambda$ is just a subgraph counting statistic based on a Poissonized random geometric graph, and enters the general framework of \cite[Section 3.4]{penrosebook}, where general $m$-dimensional CLTs are obtained for these objects. It is immediately checked that our variance estimates coincide with those stated in \cite[p. 56]{penrosebook} (for the case $k=2$), whereas our estimates in the Wasserstein distance refine the findings of \cite[Theorems 3.9 and 3.10]{penrosebook} (in the case $k=2$ and $m=1$), where no information on the rate of convergence is given. Previous references for CLTs for Poissonized disk graphs are \cite{BhGh, JJ}, where no explicit rates of convergence are provided either. A generalization of the previously described findings to general subgraph counting in a disk graph model can be found in \cite[Section 3]{LRP2}.
}
\end{rem}

\begin{figure}
\centering
\begin{tikzpicture}
 
  \draw[->,very thin] (-3.5,0) -- (3.5,0) node[above] {x};
  \draw[->,very thin] (0,-3.5) -- (0,3.5)node[right] {y};
   % \draw[<->,dashed,very thin] (0.4,-1)  --(.4,1)  node[near start, left] {1};
%      \draw[<->,very thin,dashed] (-1, 0.4) -- (-0.25,0.4) node[below] {1} -- (1,.4);
\draw  (-1,-1) [dashed ]rectangle +(2,2);

\draw[domain=1:5,smooth,variable=\x] plot (\x,{1/\x^2});
\draw[domain=1:5,smooth,variable=\x] plot (\x,{-1/\x^2});

  \draw[->,very thin] (2, 0) --  node[above=7pt]  {$f(x)$}  (2,0.25) ;

\draw[domain=-5:-1,smooth,variable=\x] plot (\x,{1/\x^2});
\draw[domain=-5:-1,smooth,variable=\x] plot (\x,{-1/\x^2});

\draw[domain=1:5,smooth,variable=\y] plot ({1/\y^2},\y);
\draw[domain=1:5,smooth,variable=\y] plot ({-1/\y^2},\y);

  \draw[->,very thin] ( 0,2) --  node[ right=9pt]  {$f(y)$}  (0.25,2) ;
  
  \draw[domain=-5:-1,smooth,variable=\y] plot ({1/\y^2},\y);
\draw[domain=-5:-1,smooth,variable=\y] plot ({-1/\y^2},\y);

%\draw[domain=-5:1,smooth,variable=\x] plot (\x,{1/\x});
%\draw[domain=-5:1,smooth,variable=\x] plot (\x,{-1/\x});

% \draw[scale=0.5,domain=1:.1:10,smooth,variable=\x] plot (\x,{1/\x}) ;

%\tikz \draw[x=1.57ex,y=1ex] (0,0) sin (1,1) sin (4,0) 
                        %    (0,1) cos (1,0) sin (2,-1) cos (3,0) sin (4,1)
;
\end{tikzpicture}

 \caption{General form of $G_{\lambda}=\{(x,y): |y|<|f(x)|\text{ and }|x|<|f(y)|\}$.}
\label{fig:G-flower}
\end{figure}
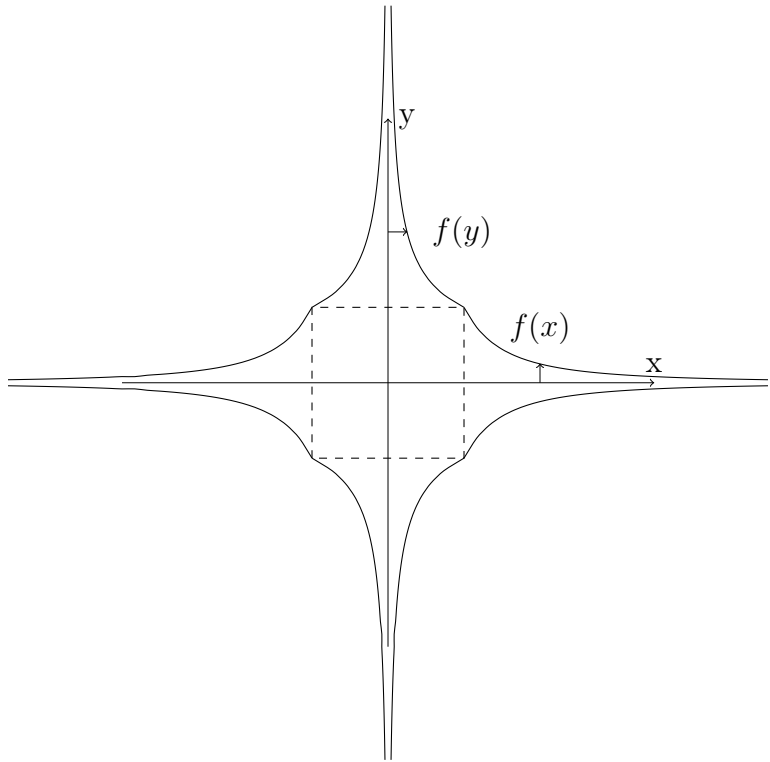

\begin{example}{\rm We present here examples of geometric stationary graphs with arbitrarily long connections and behaviors spanning the whole spectrum of possibilities indicated in Theorem \ref{thm:CLT-stat-graphs}. In all cases, condition (\ref{eq:psi-regular}) is easily checked.
We define $\fact G \subseteq \mathbb{R}^2$ as in Fig. \ref{fig:G-flower}, symmetric with respect to all axes, and consider different choices for $f$, a monotone function $\mathbb{R}_{+}^*\mapsto \mathbb{R}_{+}^*$. We set $G_{\lambda}=\alpha_{\lambda}G$ with different values for $\alpha_{\lambda}>0$. We refer the reader to Fig. \ref{fig:simuls} for results of simulations, where different values for $\lambda$ have been used. 
\begin{description}
\item[(a)] If $f(x)=1/x,\,  \alpha_{\lambda}=1$, then $\lambda\psi(\lambda)\asymp \log(\lambda)$, whence there is a CLT with an upper bound of the order of $\lambda^{-1/2}$ (Situation (i) of Theorem \ref{thm:CLT-stat-graphs}), see Fig. \ref{fig:a}.
\item[(b)] If $f(x)=1/x^2,\,  \alpha_{\lambda}=1$, then $\lambda\psi(\lambda)\asymp 1$, whence there is a CLT with an upper bound of the order of $\lambda^{-1/2}$ (Situation (ii) of Theorem \ref{thm:CLT-stat-graphs}), see Fig. \ref{fig:b}.
\item[(c)] If $f(x)=1/x^2,\, \alpha_{\lambda}=\lambda^{-1/4}$, then $\lambda\psi(\lambda)\asymp \lambda^{-1/2}$, whence there is a CLT with an upper bound of the order of $\lambda^{-1/4}$ (Situation (iii) of Theorem \ref{thm:CLT-stat-graphs}), see Fig. \ref{fig:c}.
\item[(d)] If $f(x)=1/x,\, \alpha_{\lambda}=\lambda^{-1/2}$, then $\lambda\psi(\lambda)\asymp \log(\lambda)\lambda^{-1}$, whence there is a CLT with an upper bound of the order of $\log(\lambda)^{-1}$ (Situation (iii) of Theorem \ref{thm:CLT-stat-graphs}), see Fig. \ref{fig:d}. 
\item[(e)] $f(x)=1/x^2, \alpha_{\lambda}=\lambda^{-1/2}$. According to Theorem \ref{t:poisson}, $\tilde F_{\lambda}$ converges in distribution to a multiple of a Poisson random variable, see Fig. \ref{fig:e}.
\end{description}}
\end{example}

In all cases, the convergence to a normal law goes hand in hand with the almost sure convergence of the number of connections to infinity, and with the convergence of the variance to infinity. In the case {\bf (d)}, the convergence is very slow, the number of connections behaves asymptotically like a Poisson law with parameter $\log(\lambda)$,  due to the long-range connections within the point process. In the case {\bf (e)}, the asymptotic properties of $G_\lambda$ do not yield  long range connections and the number of connections converges towards a Poisson-type limit. 

\clearpage
\newpage

\begin{subfigures} \label{fig:simuls}
 \begin{figure}\centering
\includegraphics[scale=.35]{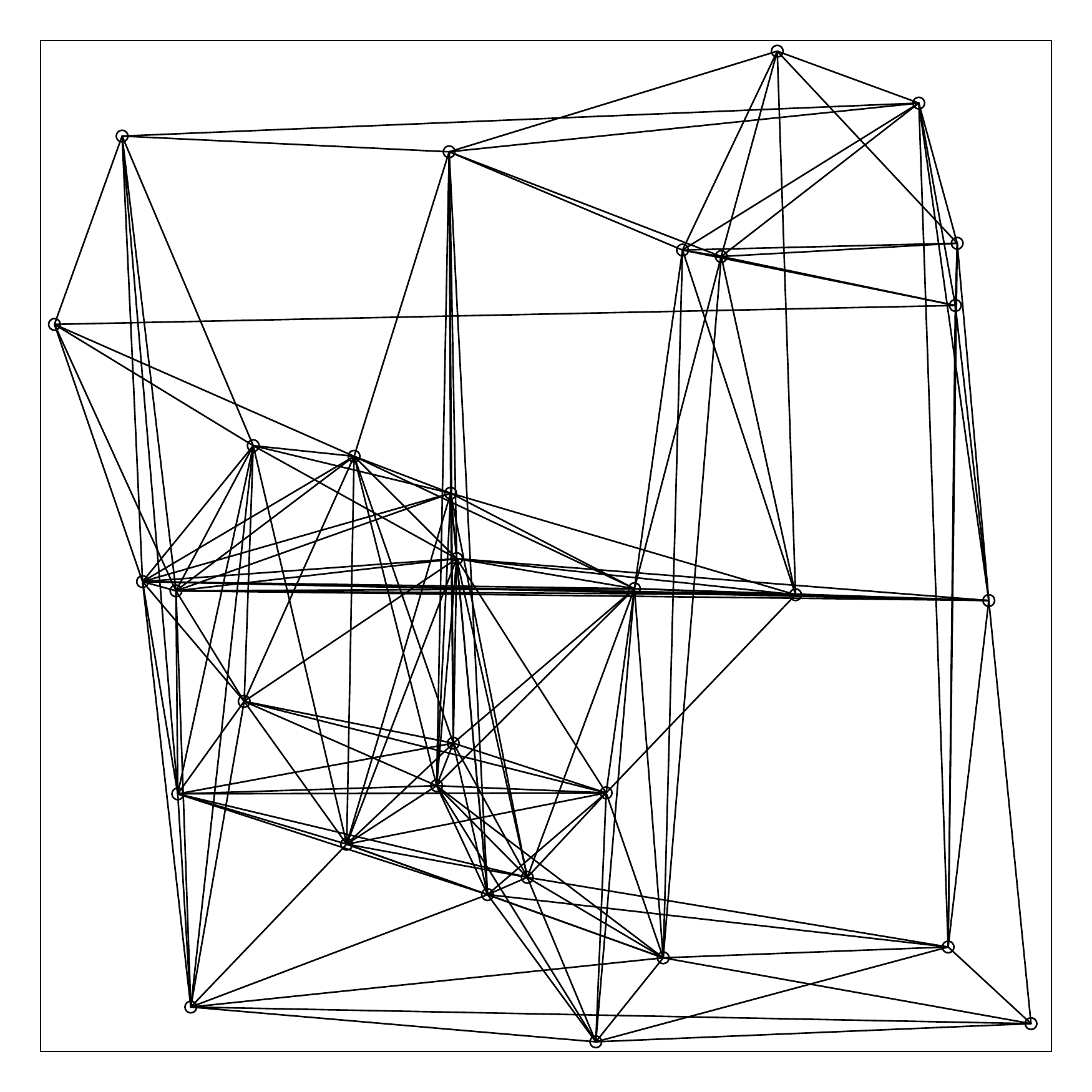}
 \includegraphics[scale=.35]{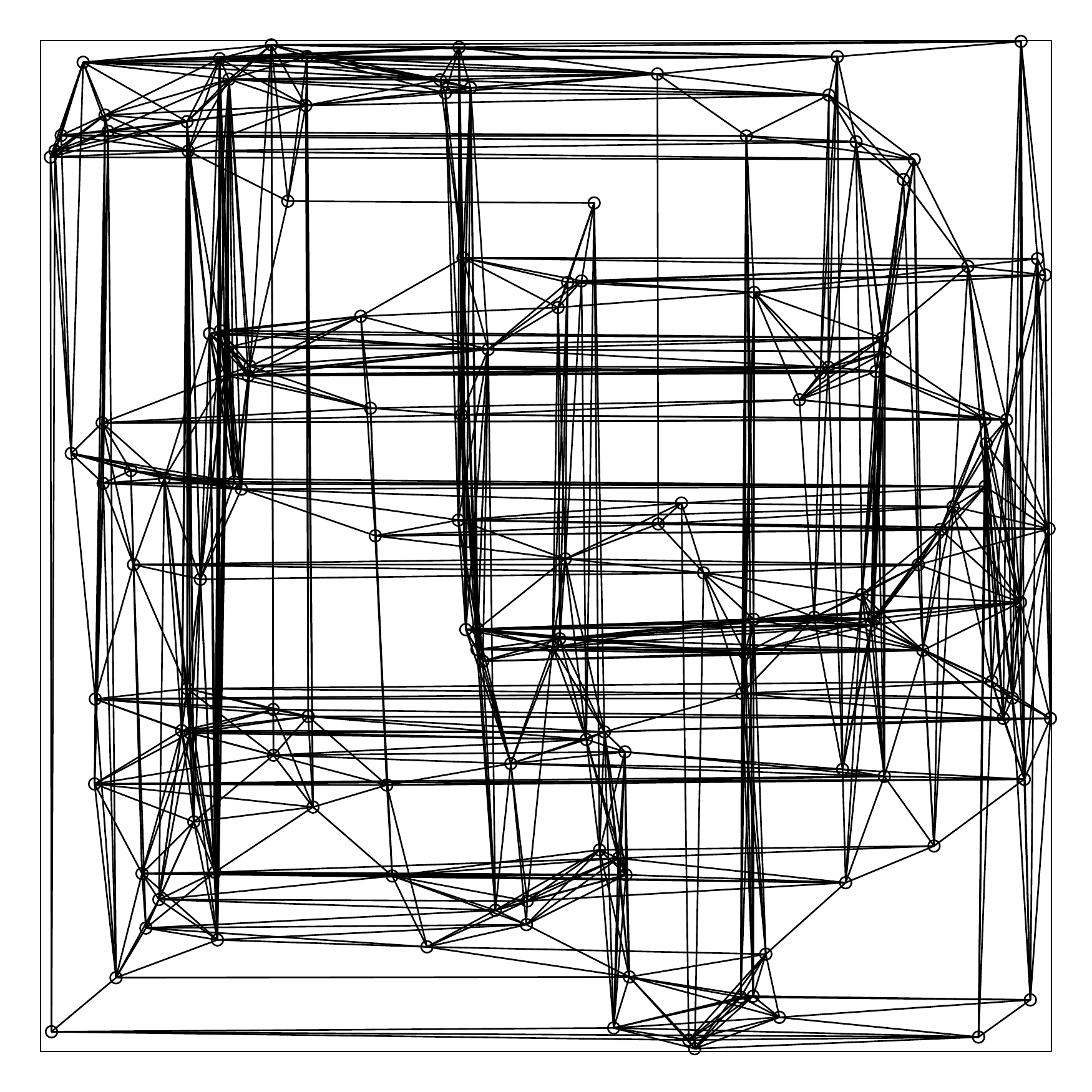}
\caption{$\lambda=25, \lambda=150$}\label{fig:a} \end{figure}
\begin{figure}\centering
\includegraphics[scale=.35]{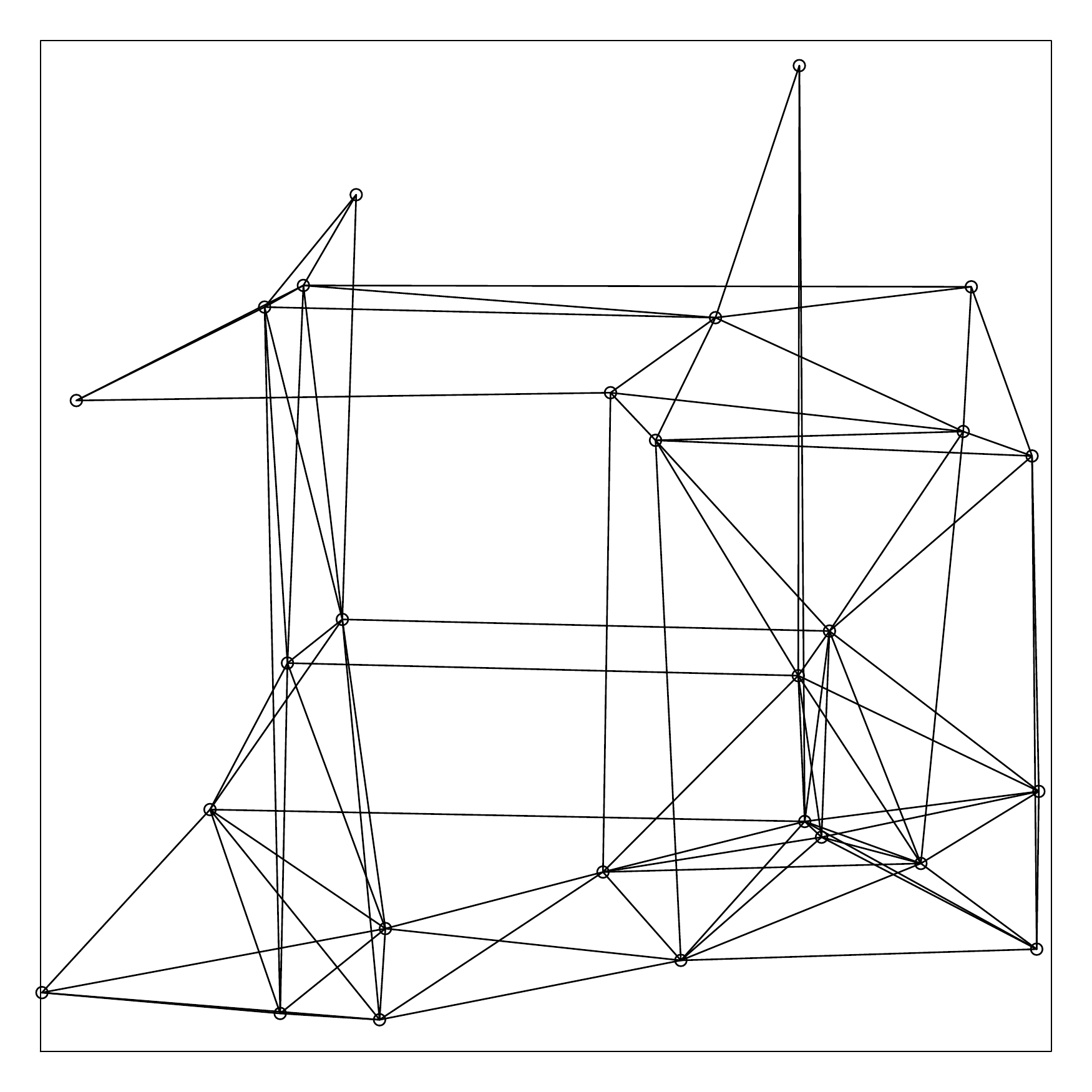}
\includegraphics[scale=.35]{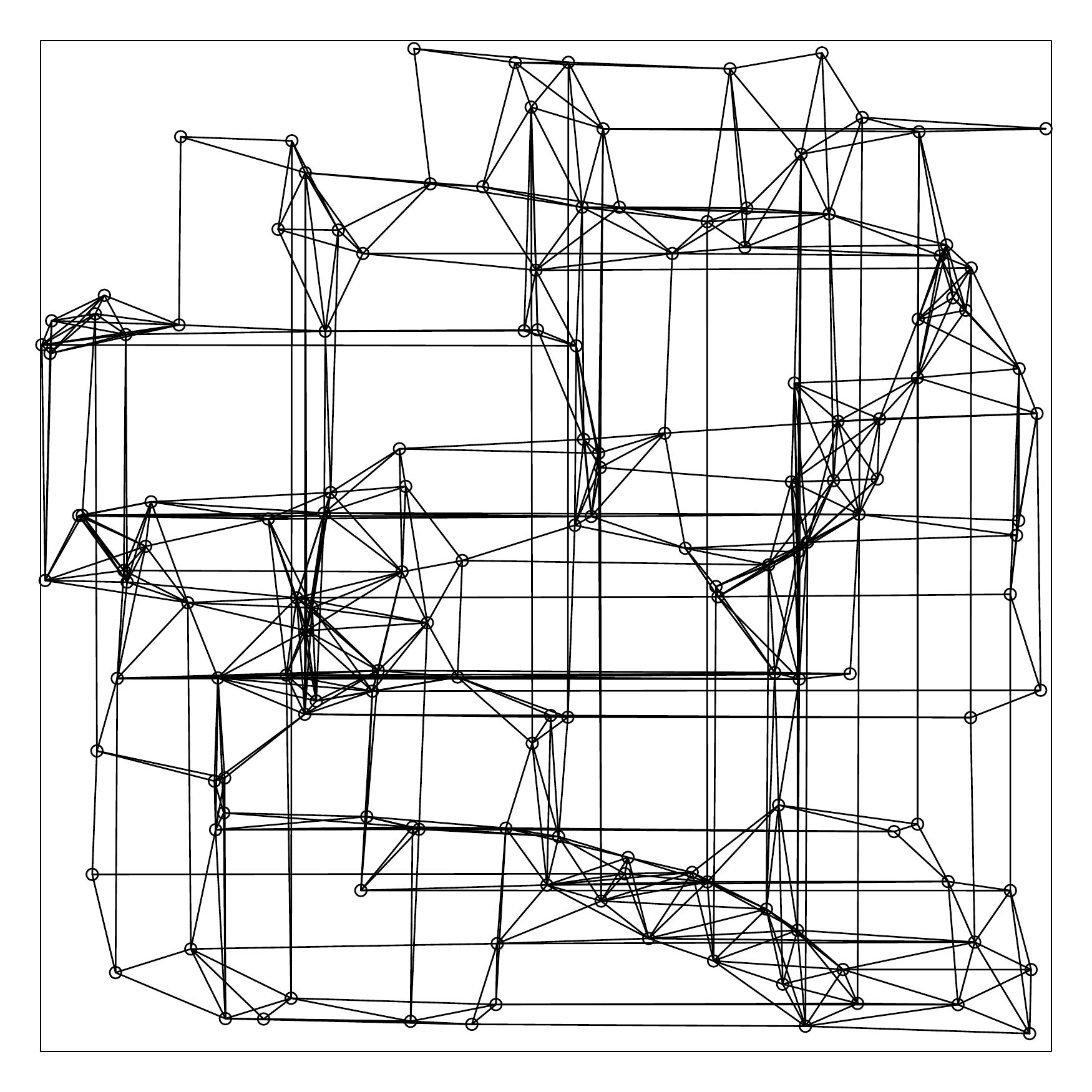}
\caption{$\lambda=25, \lambda=150$}\label{fig:b} \end{figure}
\begin{figure}\centering
\includegraphics[scale=.35]{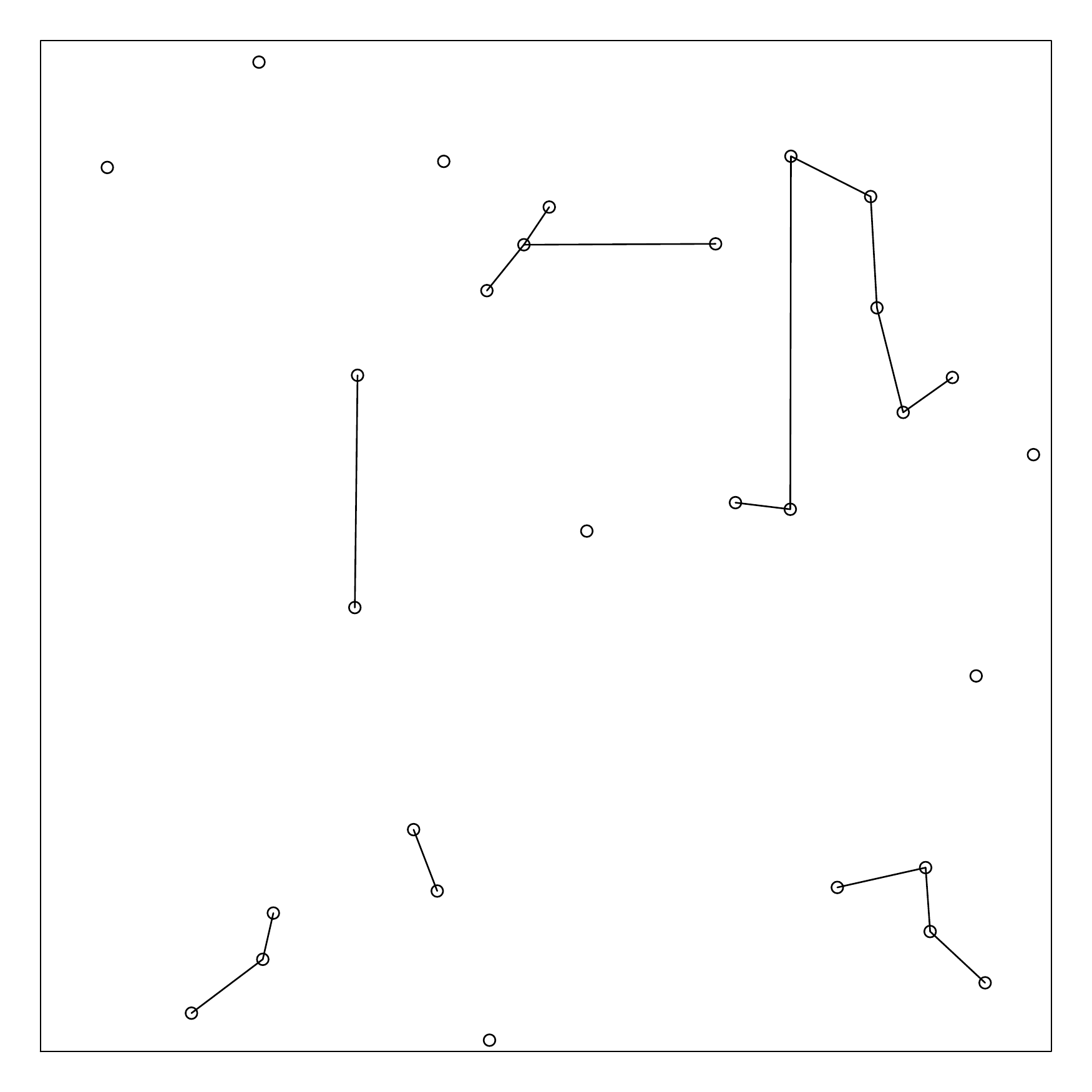}
\includegraphics[scale=.35]{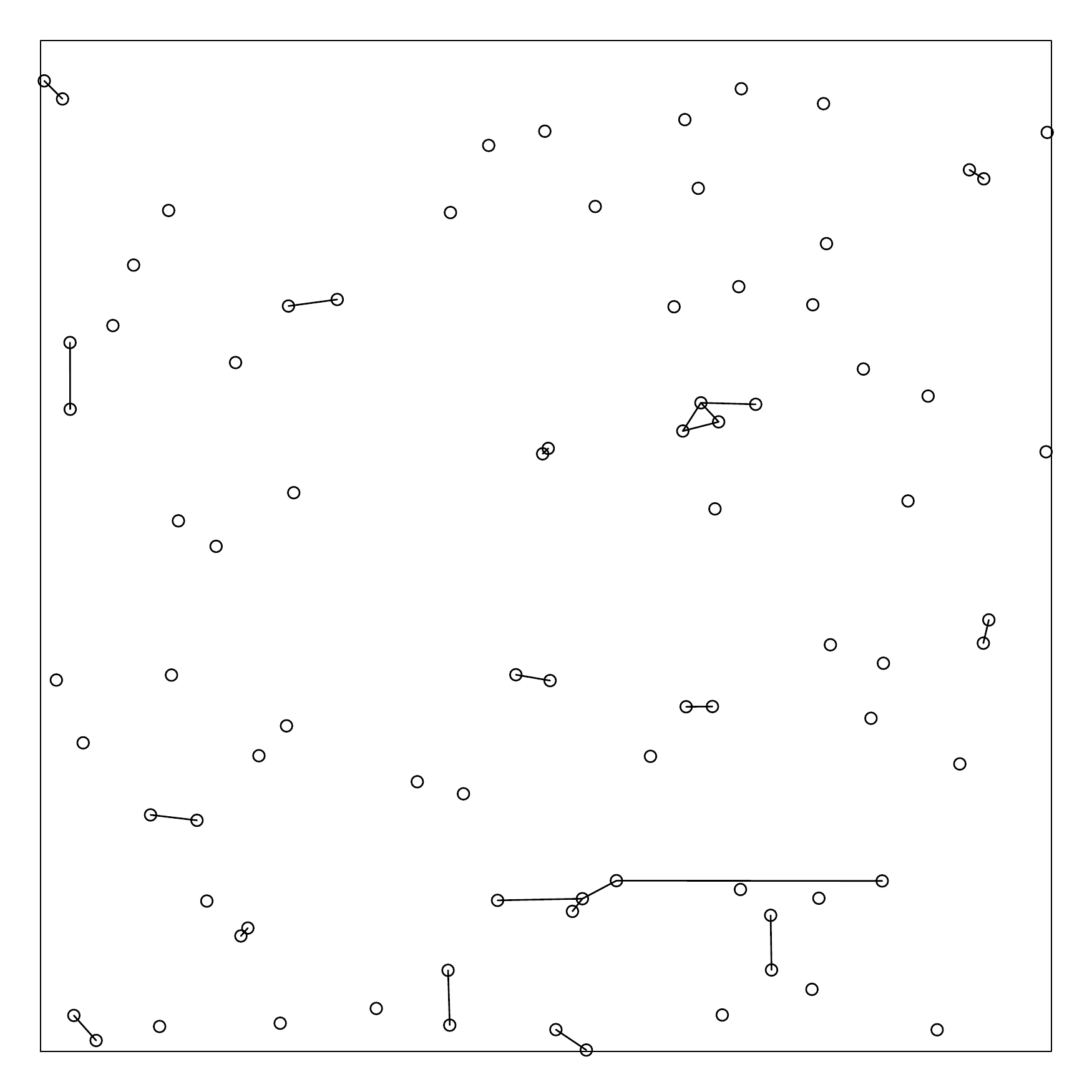}
\caption{$\lambda=25, \lambda=75$}\label{fig:c} \end{figure}
\begin{figure}\centering
\includegraphics[scale=.35]{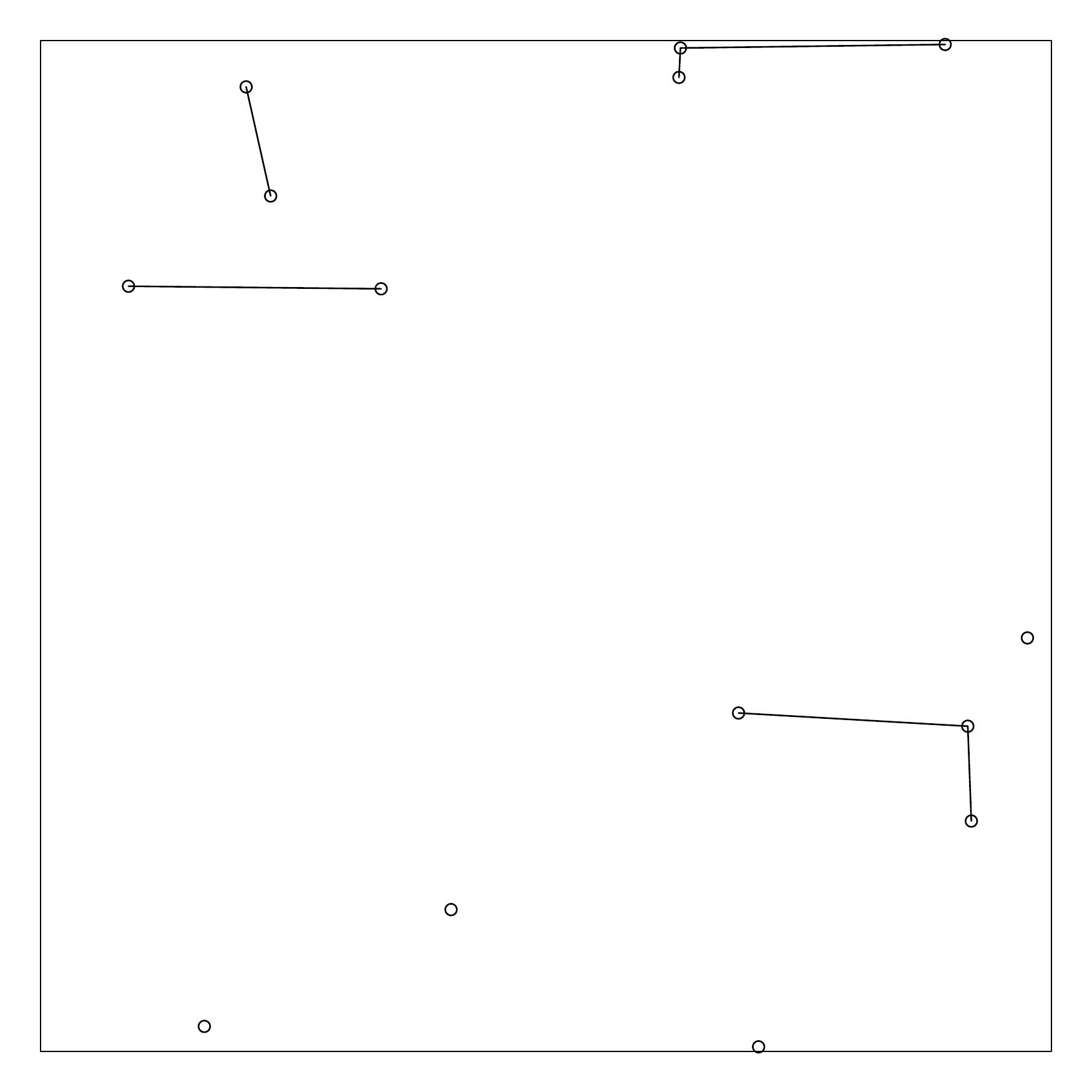}
\includegraphics[scale=.35]{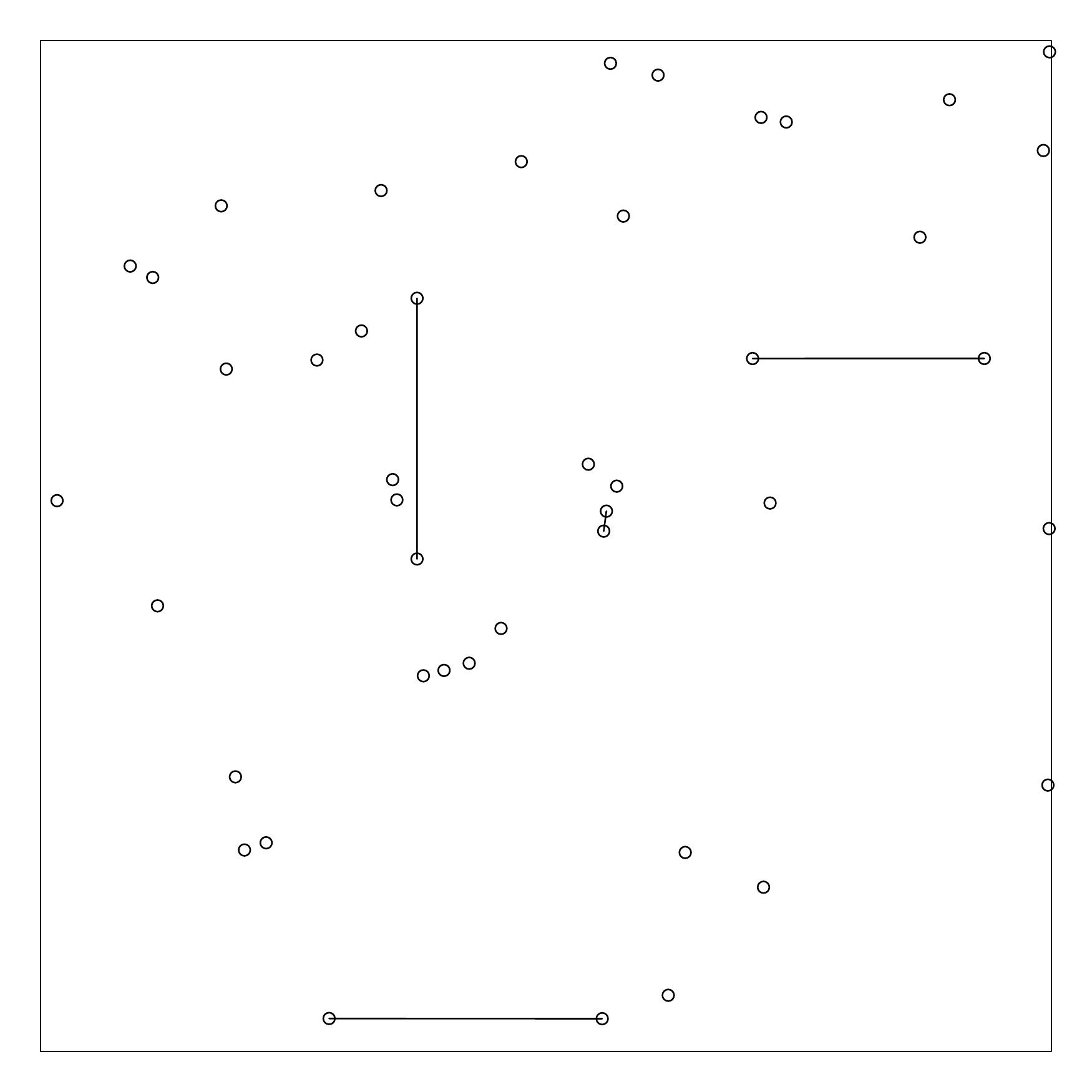}
\caption{$\lambda=10, \lambda=50$}\label{fig:d} \end{figure}
\begin{figure}\centering
\includegraphics[scale=.35]{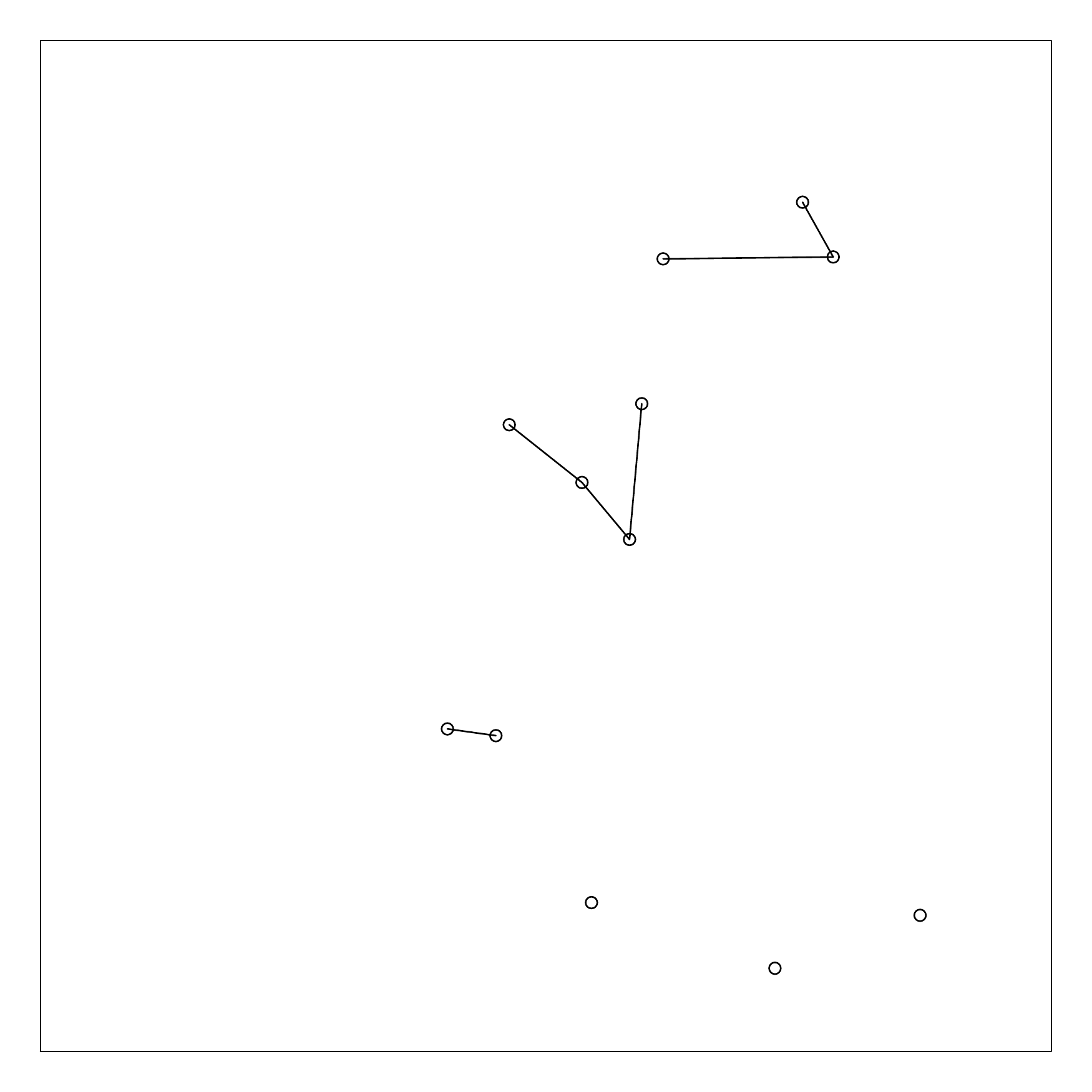}
\includegraphics[scale=.35]{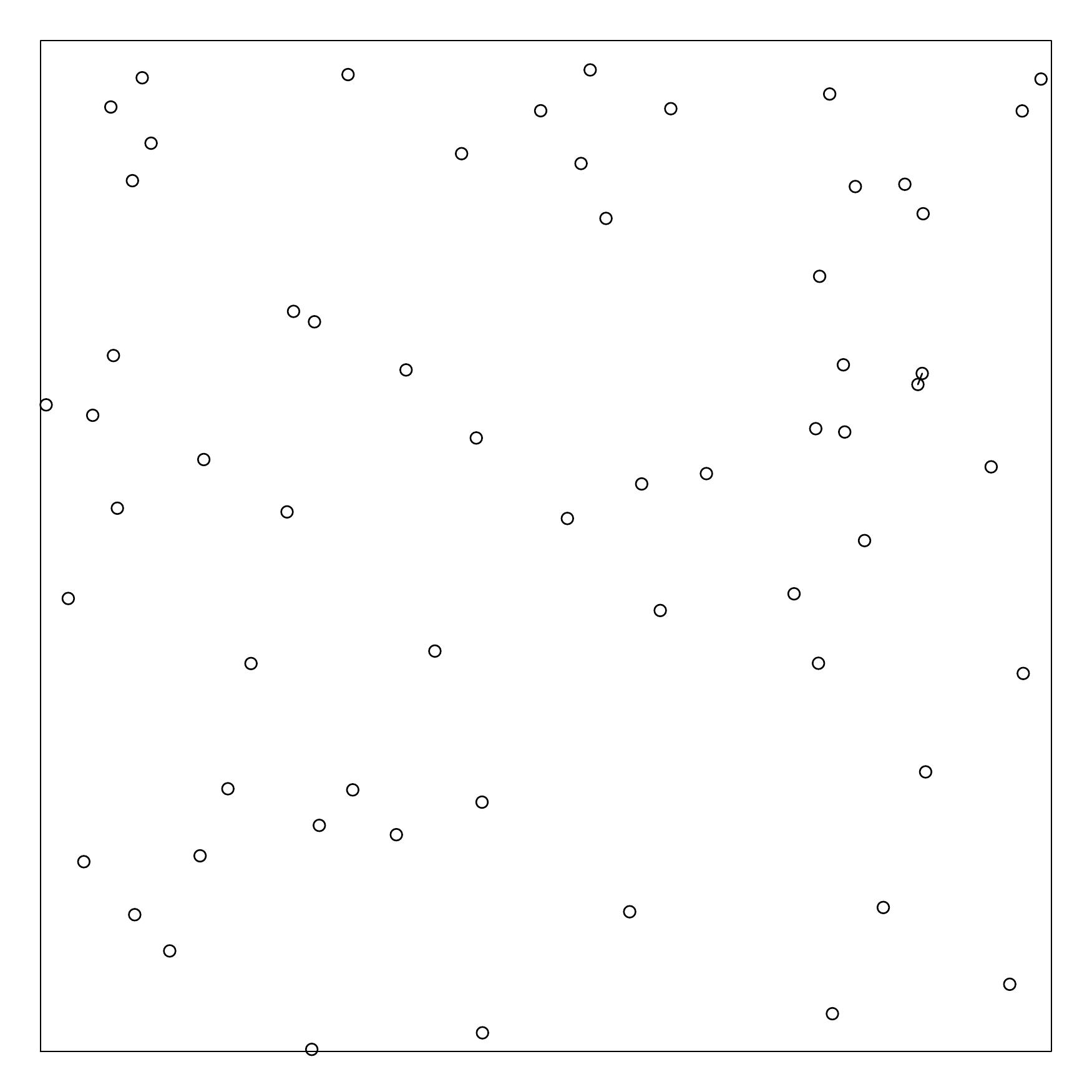}
\caption{$\lambda=10, \lambda=50$}\label{fig:e} \end{figure}
\end{subfigures}

%
%\subfloat[Subfigure 1 list of figures text][$\lambda=10, \lambda=50$]{
%\includegraphics[scale=.35]{Rplot-d-10.pdf}
%\includegraphics[scale=.35]{Rplot-d-50.pdf}
%\label{fig:d}}
%
	
%\subfloat[Subfigure 1 list of figures text][$\lambda=10, \lambda=50$]{
%\includegraphics[scale=.35]{Rplot-e-10.pdf}
%\includegraphics[scale=.35]{Rplot-e-50.pdf}
%\label{fig:e}}

%\begin{figure}[!h]
%\centering
	
%\subfloat[Subfigure 1 list of figures text][$\lambda=25, \lambda=150$]{
%\includegraphics[scale=.35]{Rplot-a-25.pdf}
%\includegraphics[scale=.35]{Rplot-a-25.pdf}
%\label{fig:a}}

%\subfloat[Subfigure 1 list of figures text][$\lambda=25, \lambda=150$]{
%\includegraphics[scale=.35]{Rplot-b-25.pdf}
%\includegraphics[scale=.35]{Rplot-b-150.pdf}
%\label{fig:b}}

%\subfloat[Subfigure 1 list of figures text][$\lambda=25, \lambda=75$]{
%\includegraphics[scale=.35]{Rplot-c-25.pdf}
%\includegraphics[scale=.35]{Rplot-c-75.pdf}
%\label{fig:c}}

%\caption{ Simulations}
%\label{fig:simuls}

%\end{figure}

%\begin{figure}[!h]
%\centering

%\subfloat[Subfigure 1 list of figures text][\textbf{( c )}: $\lambda=25, \lambda=75$]{
%\includegraphics[scale=.35]{Rplot-c-25.pdf}
%\includegraphics[scale=.35]{Rplot-c-75.pdf}
%\label{fig:c}}

%\subfloat[Subfigure 1 list of figures text][$\lambda=10, \lambda=50$]{
%\includegraphics[scale=.35]{Rplot-d-10.pdf}
%\includegraphics[scale=.35]{Rplot-d-50.pdf}
%\label{fig:d}}

%\subfloat[Subfigure 1 list of figures text][$\lambda=10, \lambda=50$]{
%\includegraphics[scale=.35]{Rplot-e-10.pdf}
%\includegraphics[scale=.35]{Rplot-e-50.pdf}
%\label{fig:e}}

%\caption{Simulations}
%\label{fig:simuls}

%\end{figure}

\bibliographystyle{plain}

\begin{section}{Appendix: Malliavin operators on the Poisson space}\label{APPENDIW}
We now define some Malliavin-type operators associated with a Poisson measure $\eta$, on the Borel space $(Z,\mathscr{Z})$, with non-atomic control measure $\mu$. We follow the work by Nualart and Vives \cite{nuaviv}.

\medskip

\noindent \underline{\bf The derivative operator $D$}. 

\smallskip

For every $F\in L^2(P)$, the derivative of $F$, $DF$ is defined as an element of $L^2(P;L^2(\mu))$, that is, of the space of the jointly measurable random functions $u:\Omega \times Z \mapsto \mathbb{R}$ such that $E \left[\int_Z u_z^2 \mu(dz) \right] <\infty$.
\begin{defi}
 \begin{enumerate}
   \item The domain of the derivative operator $D$, written  ${\rm dom} D$, is the set of all random variables $F\in L^2(P)$ admitting a chaotic decomposition (\ref{chaos}) such that
$$ \sum_{k\geq 1} k k!\|f_k \|^2_{L^2(\mu^k)} < \infty ,$$
   \item For any $F\in {\rm dom}D$, the random function $z \mapsto D_z F$ is defined by
$$ D_z F= \sum_{k \geq 1}^{\infty} k I_{k-1}(f_k(z,\cdot)) .$$
 \end{enumerate}
\end{defi}

\medskip

\noindent \underline{\bf The divergence operator $\delta$}.  

\smallskip

Thanks to the chaotic representation property of $\eta$, every random function
$u \in L^2(P,L^2(\mu))$ admits a unique representation of the type
\begin{equation} \label{skor}
 u_z = \sum_{k \geq 0}^{\infty}  I_{k}(f_k(z,\cdot)) ,\,\, z\in Z,
\end{equation}
where the kernel $f_k$ is a function of $k+1$ variables, and $f_k(z,\cdot)$ is an element of $L^2_s(\mu^k)$. The {\sl divergence operator} $\delta(u)$ maps a random function $u$ in its domain to an element of $L^2(P)$.\\

\begin{defi}
\begin{enumerate}
  \item  The domain of the divergence operator, denoted by  ${\rm dom} \delta$, is the collection of all $u\in L^2(P,L^2(\mu))$ having the above chaotic expansion (\ref{skor}) satisfied the condition:
$$ \sum_{k\geq 0}  (k+1)! \|f_k \|^2_{L^2(\mu^(k+1))} < \infty. $$
  \item For $u\in {\rm dom}\delta$, the random variable $\delta(u)$ is given by
      $$ \delta (u) = \sum_{k\geq 0} I_{k+1}(\tilde{f}_k), $$
      where $\tilde{f}_k$ is the canonical symmetrization of the $k+1$ variables function $f_k$.
\end{enumerate}
\end{defi}
As made clear in the following statement, the operator $\delta$ is indeed the adjoint operator of $D$.
\begin{lemme}[Integration by parts]\label{L : IBP}
 For every $G\in {\rm dom} D$ and $u\in {\rm dom} \delta$, one has that
$$ E[G \delta(u)] = E[\langle D G, u \rangle_{L^2(\mu)}]. $$
\end{lemme}
The proof of Lemma \ref{L : IBP} is detailed e.g. in \cite{nuaviv}.\\

\medskip

\noindent \underline{\bf The Ornstein-Uhlenbeck generator $L$}.

\smallskip

\begin{defi}
\begin{enumerate}
  \item  The domain of the Ornstein-Uhlenbeck generator, denoted by  ${\rm dom} L$, is the collection of all $F \in L^2(P)$ whose chaotic representation \label{chao} verifies the condition:
$$ \sum_{k\geq 1}  k^2 k! \|f_k \|^2_{L^2(\mu^k)} < \infty $$
  \item The Ornstein-Uhlenbeck generator $L$ acts on random variable $F\in {\rm dom}L$ as follows:
      $$ LF = - \sum_{k\geq 1} k I_{k}(f_k) .$$
\end{enumerate}
\end{defi}

\medskip

\noindent \underline{\bf The pseudo-inverse of $L$}.

\smallskip

\begin{defi}
\begin{enumerate}
  \item  The domain of the pseudo-inverse of the Ornstein-Uhlenbeck generator, denoted by $L^{-1}$, is the space $L^2_0(P)$ of \it{centered} random variables in $L^2(P)$.
  \item For $F = \sum\limits_{k\geq 1} I_k (f_k) \in L^2_0(P)$ , we set
      $$ L^{-1}F = - \sum_{k\geq 1} \cfrac{1}{k} I_{k}(f_k). $$
\end{enumerate}
\end{defi}
\end{section}

\end{document}